\documentclass[12pt]{amsart}

\usepackage{amsmath}
\usepackage{amsthm}
\usepackage{amssymb}
\usepackage{eucal}
\usepackage{times}
\usepackage{euler}
\usepackage{eufrak}

\usepackage[all]{xypic}

\usepackage[dvips]{graphicx}

\author{Piotr \'Sniady}

\address{Institute of Mathematics \\
University of Wroclaw \\ pl.~Grunwaldzki~2/4 \\ 50-384 Wroclaw
\\ Poland}

\email{Piotr.Sniady@math.uni.wroc.pl}

\title[Gaussian fluctuations of characters of symmetric groups]%
{Gaussian fluctuations of characters of symmetric groups and of
Young diagrams}

\numberwithin{equation}{section} \numberwithin{figure}{section}

\theoremstyle{plain}
\newtheorem{lemma}{Lemma}[section]
\newtheorem{theorem}[lemma]{Theorem}
\newtheorem{theoremanddefinition}[lemma]{Theorem and Definition}
\newtheorem{proposition}[lemma]{Proposition}
\newtheorem{corollary}[lemma]{Corollary}

\theoremstyle{definition}

\theoremstyle{remark}
\newtheorem*{remark}{Remark}
\newtheorem{example}[lemma]{Example}

\newcommand{\A}{{\mathfrak{A}}}

\newcommand{\E}{{\mathbb{E}}}
\newcommand{\C}{{\mathbb{C}}}
\newcommand{\R}{{\mathbb{R}}}

\newcommand{\M}{{\mathcal{M}}}

\newcommand{\Sn}[1]{{{S}}_{#1}}

\newcommand{\cytat}[1]{}

\newcommand{\podpotega}[1]{^{\underline{#1}}}

 \DeclareMathOperator{\Cov}{Cov}

 \DeclareMathOperator{\fat}{fat}

 \DeclareMathOperator{\Tr}{Tr}

 \DeclareMathOperator{\genus}{genus}

\begin{document}

\begin{abstract}
We study asymptotics of reducible representations of the symmetric
groups $S_q$ for large $q$. We decompose such a representation as a
sum of irreducible components (or, alternatively, Young diagrams)
and we ask what is the character of a randomly chosen component (or,
what is the shape of a randomly chosen Young diagram). Our main
result is that for a large class of representations the fluctuations
of characters (and fluctuations of the shape of the Young diagrams)
are asymptotically Gaussian; in this way we generalize Kerov's
central limit theorem. The considered class consists of
representations for which the characters almost factorize and this
class includes, for example, left-regular representation (Plancherel
measure), tensor representations. This class is also closed under
induction, restriction, outer product and tensor product of
representations. Our main tool in the proof is the method of genus
expansion, well known from the random matrix theory.
\end{abstract}

\maketitle

\section{Introduction}

\subsection{Representations of large symmetric groups}

Irreducible representations of the symmetric groups $S_q$ are
indexed by Young diagrams and nearly all questions about them, such
as values of the characters or decomposition into irreducible
components of a restriction, induction, tensor product or outer
product of representations can be answered by combinatorial
algorithms such as Murnaghan­-Nakayama formula or the
Littlewood­-Richardson rule. Unfortunately, these exact
combinatorial tools become very complicated and cumbersome when the
size of the symmetric group $S_q$ tends to infinity. For example, a
restriction of an irreducible representation consists typically of a
very large number of Young diagrams and listing them all does not
give much insight into their structure. In order to deal with such
questions in the asymptotic region when $q\to\infty$ we should be
more modest and ask questions of a more statistical flavor: what is
the typical shape of a Young diagram contributing to a given
representation? what are the fluctuations of the Young diagrams
around the most probable shape?

In this article we are interested in the situation when---speaking
informally---a typical Young diagram contributing to the considered
representation of $S_q$ has at most $O(\sqrt{q})$ rows and columns.
The first results in this direction concerned the left-regular
representation (or equivalently, Plancherel measure on Young
diagrams): Vershik and Kerov \cite{VershikKerov1977} and Logan and
Shepp \cite{LoganShepp} found the shape of a typical Young diagram
which contributes to the left-regular representation and later Kerov
\cite{Kerov1993gaussian} announced that fluctuations of Young
diagrams contributing to the left-regular representation around
their limit shape are Gaussian (for the complete proof together with
a detailed history of the result we refer to the article by Ivanov
and Olshanski \cite{IvanovOlshanski2002}). A non-commutative version
of this result was given by Hora
\cite{Hora2002noncommutativeaspect,Hora2003noncommutativeKerov}.
Biane \cite{Biane1998,Biane2001approximate,BianeCharacters}
considered a more general case, namely representations with
approximate factorization of characters and proved that the shape of
the typical Young diagram contributing to such representations can
be described by the means of Voiculescu's free probability theory
\cite{VoiculescuDykemaNica}.

In this article we present a condition on factorization of
characters which is strong enough to ensure Gaussian fluctuations of
Young diagrams and of characters and which is weak enough to be very
common among naturally arising representations. In this way we prove
a generalization of Kerov's central limit theorem for a very wide
class of representations.

\subsection{Genus expansion, random matrices and free probability}
Free probability of Voiculescu \cite{VoiculescuDykemaNica} is a
non-commutative probability theory which turned out to be very
successful in describing random matrices. The combinatorial
structure behind this theory is the lattice of non-crossing
partitions \cite{Kreweras} and the corresponding notion of free
cumulants \cite{Speicher1997}. It was Biane \cite{Biane1998} who
realized that the same structure describes the leading terms in the
asymptotic description of representations of symmetric groups. In
this way some notions concerning random matrices were matched (via
free probability theory) to some notions concerning representations
of symmetric groups.

However, the connection between the random matrix theory and the
representation of symmetric groups is much deeper than just the
connection to free probability and non-crossing partitions. Example
of such a connection was given by Okounkov
\cite{Okounkov2000randompermutations} who showed that the joint
distribution of the largest eigenvalues of a GUE random matrix
coincides (after appropriate scaling) with the joint distribution of
the longest rows of a Young diagram distributed according to the
Plancherel measure. For proving such results it is not enough to
consider the first-order approximation given by non-crossing
partitions and one has to use exact formulas. Such formulas for the
moments of large classes of random matrices were known for a long
time and they can be viewed as series indexed by two-dimensional
surfaces \cite{Zvonkin1995matrixintegrals} and Okounkov
\cite{Okounkov2000randompermutations} found their counterpart for
random Young diagrams distributed according to the Plancherel
measure. The asymptotic behavior of a term in this expansion depends
only on the topology of the surface and for this reason such
formulas are called genus expansions. It became clear that the
origin of the similarities between the random matrix theory and the
theory of the representations of the symmetric groups is the common
structure of the genus expansion.

In our recent work \cite{Sniady2003pushing} we pointed out that the
genus expansion method can be applied not only to the Plancherel
measure but to a wide class of representations. This method is also
the main tool in the proofs in this article.

\subsection{Higher-order free probability}
\label{subsec:higherorder} As we mentioned above, the first order
approximation is not sufficient to calculate the fluctuations of
characters and of Young diagrams and therefore such fluctuations
cannot be described in the framework of the (usual) free probability
theory. On the other hand, in a series of articles
\cite{HigherOrderFreeness1,HigherOrderFreeness2,HigherOrderFreeness3}
it was demonstrated that by considering some more complicated
versions of non-commutative probability spaces it is possible to
describe fluctuations of random matrices in the framework of, so
called, higher order free probability. In this theory the
non-crossing partitions are replaced by a more general object,
namely annular non-crossing partitions.  Results presented in this
article suggest that it should be possible to describe in this
framework also the fluctuations of Young diagrams and we will deal
with this problem in a forthcoming article.

\subsection{Factorization of characters}
Biane \cite{Biane2001approximate} proved (under some mild technical
assumptions) that Young diagrams contributing to some
finite-dimensional reducible representation of the symmetric group
$S_q$ will concentrate around some limit shape if and only if the
normalized character of the representation
$$\chi(\pi):=\frac{\Tr \rho(\pi)}{\Tr \rho(e)}$$
approximately factorizes, i.e.\ informally speaking
\begin{equation}
\label{eq:approximate} \chi( \sigma_1 \cdots \sigma_n) \approx
\chi(\sigma_1) \cdots \chi(\sigma_n)
\end{equation}
for all permutations $\sigma_1,\dots,\sigma_n$ with disjoint
supports (to be more precise: we consider a sequence of
representations $(\rho_q)$ where $\rho_q$ is a representation of
$S_q$ and the approximate equality \eqref{eq:approximate} should
hold in the limit $q\to\infty$). The result of Biane can be viewed
as an analogue of the law of large numbers while the results
presented in this article are an analogue of the central limit
theorem; one can ask therefore which condition should replace
\eqref{eq:approximate} in order to prove such stronger results.

We will not be very far from the truth when we say that for the
results of Biane \cite{Biane2001approximate} it is enough to assume
some version of \eqref{eq:approximate} for $n=2$. If we treat
permutations $\sigma_1,\sigma_2$ as random variables and the
normalized character $\chi$ as an expectation $\E$ then this
condition can be equivalently written as a condition for the
covariance:
$$ \Cov(\sigma_1,\sigma_2)=\E( \sigma_1
\sigma_2)-\E(\sigma_1) \E(\sigma_2) \approx 0.$$ Covariance is a
special case of a more general probabilistic notion of a cumulant
(we will recall the necessary definitions in Section
\ref{subsec:cumulants}) therefore it is quite natural to expect that
the correct condition for the factorization of characters should
involve all cumulants and indeed in this article we prove that the
following condition is sufficient for our purposes (Theorem and
Definition \ref{theo:main}): for any cycles
$\sigma_1,\dots,\sigma_n$ with disjoint supports we assume that
\begin{equation}
\label{eq:approximate2} k_n(\sigma_1,\dots,\sigma_n) = O\big(q^{-
\frac{|\sigma_1|+\cdots+|\sigma_n|+2(n-1)}{2} } \big),
\end{equation}
where $|\pi|$ denotes the minimal number of factors needed to write
the permutation $\pi$ as a product of transpositions.

It should be stressed that the decay of the cumulants in the
condition \eqref{eq:approximate2} carries a strong resemblance to
the decay of the cumulants of the entries of many interesting
classes of random matrices \cite{Collins2002} and, as we mentioned
in Section \ref{subsec:higherorder}, it is not an accident.

\subsection{Fluctuations of Young diagrams}

Every finite-dimensional representation $\rho_q$ of the symmetric
group $S_q$ defines a canonical probability measure on the Young
diagrams contributing to $\rho_q$ given as follows: probability of a
Young diagram $\lambda$ should be proportional to the total
dimension of all irreducible components of type $[\lambda]$. Our
goal is to consider some interesting function $f$ on the set of the
Young diagrams with $q$ boxes and to study the distribution of the
random variable $f(\lambda)$. In principle, the information about
the characters such as \eqref{eq:approximate2} should be sufficient
to compute the distribution of random variables $f(\lambda)$ for
reasonable functions $f$, however this relation is not very direct.
An analogue of this situation can be found in the random matrix
theory, where the knowledge of the joint distribution of the entries
of a random matrix should be enough to find the joint distribution
of the eigenvalues, however the actual calculation might be quite
involved.

In this article we will show that the joint distribution of the
random variables of the form $f(\lambda)$ converges to a Gaussian
one if $f$ is the value of the corresponding irreducible character
on a prescribed permutation or some functional describing the shape
of $\lambda$.


\subsection{Overview of this article}

In Section \ref{sec:preliminaries} we present briefly all necessary
notions needed to state the main result and its applications. In
Section \ref{sec:representations} we present the main result:
Theorem and Definition \ref{theo:main} where four equivalent
conditions are given which ensure Gaussian fluctuations of the
characters and of the shape of the Young diagrams. We also show that
the considered class of representations has many interesting
examples and that it is closed under some natural operations such as
induction, restriction, outer product and tensor product of
representation. Section \ref{sec:Proof-main-result} contains the
proof of the main result. Finally, in Section
\ref{sec:Proof-Theorem-?????} we proof some technical results used
in the proof of the main theorem.


\section{Preliminaries}
\label{sec:preliminaries}

\subsection{Normalized conjugacy class indicators}

\label{subsec:definicjasigma}

Let integer numbers $k_1,\dots,k_m\geq 1$ be given. We define the
normalized conjugacy class indicator to be a central element in the
group algebra $\C(S_q)$ given by
\cite{KerovOlshanski1994,BianeCharacters,Sniady2003pushing}
\begin{equation}
\label{eq:definicjasigma} \Sigma_{k_1,\dots,k_m}=
\sum_{a} (a_{1,1},a_{1,2},\dots,a_{1,{k_1}}) \cdots
(a_{m,1},a_{m,2},\dots,a_{m,k_m}), \end{equation}
where the sum runs over all one--to--one functions
$$a:\big\{ \{r,s\}: 1\leq r\leq m, 1\leq s\leq {k_r}\big\}\rightarrow \{1,\dots,q\}$$ and
$(a_{1,1},a_{1,2},\dots,a_{1,{k_1}}) \cdots
(a_{m,1},a_{m,2},\dots,a_{m,k_m})$ denotes the product of disjoint
cycles. Of course, if $q<k_1+\cdots+k_m$ then the above sum runs
over the empty set and $\Sigma_{k_1,\dots,k_m}=0$.

In other words, we consider a Young diagram with the rows of the
lengths $k_1,\dots,k_m$ and all ways of filling it with the elements
of the set $\{1,\dots,q\}$ in such a way that no element appears
more than once. Each such a filling can be interpreted as a
permutation when we treat the  rows of the Young tableau as disjoint
cycles.

It follows that $\Sigma_{k_1,\dots,k_m}$ is a linear combination of
permutations which in the cycle decomposition have cycles of length
$k_1,\dots,k_m$ (and, additionally, $q-(k_1+\cdots+k_m)$
fix-points). Each such a permutation appears in
$\Sigma_{k_1,\dots,k_m}$ with some positive integer multiplicity
depending on the symmetry of the tuple $k_1,\dots,k_m$.

\subsection{Disjoint product}
\label{subsec:disjoint} Usually by a product of normalized conjugacy
classes $\Sigma_{k_1,\dots,k_m}$ we understand their product as
elements of the group algebra $\C(S_q)$. However, sometimes it is
convenient to consider their disjoint product defined by
\begin{equation}
\label{eq:disjoint} \Sigma_{k_1,\dots,k_m} \bullet
\Sigma_{l_1,\dots,l_n}= \Sigma_{k_1,\dots,k_m,l_1,\dots,l_n}.
\end{equation}

\begin{remark}
The readers familiar with the notion of partial permutations of
Ivanov and Kerov \cite{IvanovKerov1999} will see that
\eqref{eq:disjoint} is compatible with the following definition of
the product $\alpha_1 \bullet \alpha_2$, when $\alpha_1,\alpha_2$
are partial permutations: \begin{equation}
\label{eq:definition-of-disjoint} \alpha_1\bullet \alpha_2 =
\begin{cases} \alpha_1 \alpha_2 & \text{if supports of $\alpha_1$
and $\alpha_2$ are disjoint,} \\ 0 & \text{otherwise}. \end{cases}
\end{equation} Further discussion can be found in Section
\ref{subsec:algebraofconjugacy}.
\end{remark}

\subsection{Classical cumulants}
\label{subsec:cumulants} The notion of cumulants (or
semi-invariants) was introduced to describe the convolution of
measures in the classical probability theory. For more details about
cumulants in the classical and non-commutative probability theory we
refer to overview articles
\cite{Hald-EarlyHistoryCumulants,MattnerWhat-are-cumulants,Lehner2002cumulantsI}.

If $X_1,\dots,X_n$ are random variables we define their (classical)
cumulant to be an appropriate coefficient of the formal expansion of
the logarithm of the multidimensional Fourier transform:
$$ k(X_1,\dots,X_n)= \frac{\partial^n}{\partial t_1 \cdots
\partial t_n}\bigg|_{t_1=\cdots=t_n=0} \log \E e^{t_1 X_1+\cdots+t_n
X_n}. $$ Cumulant is linear with respect to each of its arguments.
If random variables $X_1,\dots,X_n$ can be split into two groups
such that the random variables from the first class are independent
with the random variables from the second class then their cumulant
vanishes: $k(X_1,\dots,X_n)=0$.

The first two cumulants \begin{align*} k(X)& =\E X, \\ k(X_1,X_2)&
=\E (X_1 X_2)- \E (X_1) \E(X_2)= \Cov(X_1,X_2)
\end{align*}
coincide with the mean value and the covariance.

\subsection{Elements of the group algebra as random variables}
Let us fix some finite-dimensional representation $\rho_q$ of the
symmetric group $S_q$. We can treat any commuting family of elements
of the group algebra $\C(S_q)$ as a family of random variables
equipped with the mean value given by the normalized character:
\begin{equation}
\label{eq:wartoscoczekiwana} \E X:= \chi^{\rho_q}(X)= \frac{\Tr
\rho_q(X)}{\Tr \rho_q(1)}.
\end{equation}
It should be stressed that in the general case we treat elements of
$\C(S_q)$ as random variables only on a purely formal level; in
particular we do not treat them as functions on some Kolmogorov
probability space.

Usually by the product of such random variables we understand the
natural product in the group algebra $\C(S_q)$ and we denote the
resulting cumulants (called natural cumulants) by
$k(X_1,\dots,X_n)$. However, sometimes it is more convenient to take
as the product of random variables the disjoint product $\bullet$;
we denote the resulting cumulants (called disjoint cumulants) by
$k^{\bullet}(X_1,\dots,X_n)$.

\subsection{Canonical probability measure on Young diagrams associated to a representation}
There is a special case when it is possible to give a truly
probabilistic interpretation to \eqref{eq:wartoscoczekiwana}: it is
when for the family of random variables we take the center of
$\C(S_q)$ and the multiplication is the natural product in the group
algebra $\C(S_q)$. The center of $\C(S_q)$ is isomorphic (via
Fourier transform) to the algebra of functions on Young diagrams and
the expected value \eqref{eq:wartoscoczekiwana} corresponds under
this isomorphism to the probability measure on Young diagrams with
$q$ boxes such that the probability of $\lambda$ is proportional to
the total dimension of the irreducible components of type
$[\lambda]$ in $\rho_q$.

\subsection{Generalized Young diagrams} \index{generalized Young
diagram} \index{Young diagram!generalized}
%

\begin{figure}[tb]
\includegraphics{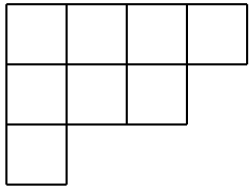}
\caption[Example of a Young diagram]{Young diagram associated with a
partition $8=4+3+1$.} \label{fig:young1}
\end{figure}

\begin{figure}[tb]
\includegraphics{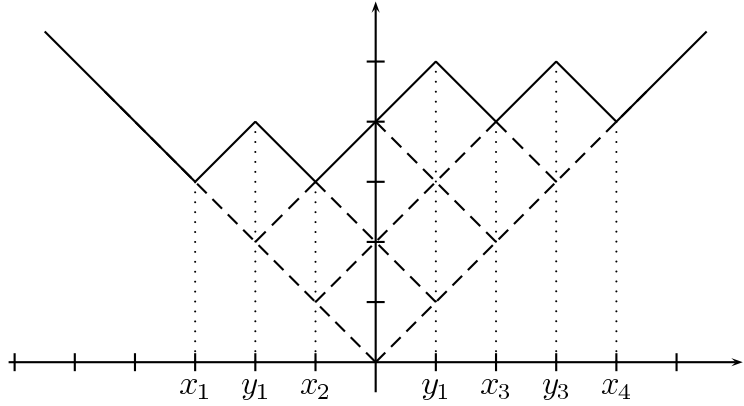}
\caption[Example of a generalized Young diagram]{Generalized Young
diagram associated with a partition $8=4+3+1$.} \label{fig:young2}
\end{figure}

Let $\lambda$ be a Young diagram. We assign to it a piecewise affine
function $\omega^\lambda:\R\rightarrow\R$ with slopes $\pm 1$, such
that $\omega^\lambda(x)=|x|$ for large $|x|$ as it can be seen on
the example from Figure \ref{fig:young2}. By comparing Figure
\ref{fig:young1} and Figure \ref{fig:young2} one can easily see that
the graph of $\omega^\lambda$ can be obtained from the graphical
representation of the Young diagram by an appropriate mirror image,
rotation and scaling by the factor $\sqrt{2}$. We call
$\omega^\lambda$ the generalized Young diagram associated with the
Young diagram $\lambda$ \cite{Kerov1993transition,
Kerov1998interlacing,Kerov1999differential}.


The class of generalized Young diagrams consists of all functions
$\omega:\R\rightarrow\R$ which are Lipschitz with constant $1$ and
such that $\omega(x)=|x|$ for large $|x|$ and of course not every
generalized Young diagram can be obtained by the above construction
from some Young diagram $\lambda$.

The setup of generalized Young diagrams is very useful in the study
of the asymptotic properties since it allows us to define easily
various notions of convergence of the Young diagram shapes.

\subsection{Functionals of the shape of Young diagrams}
The main result of this article is that the fluctuations of the
shape of some random Young diagrams converge (after some rescaling)
to a Gaussian distribution. Since the space of (generalized) Young
diagrams is infinite-dimensional therefore we need to be very
cautious when dealing with such statements. In fact, we will
consider a family of functionals on Young diagrams and we show that
the joint distribution of each finite set of these functionals
converges to the Gaussian distribution.

The functionals mentioned above are given as follows: for a Young
diagram $\lambda$ and the corresponding generalized Young diagram
$\omega$ we denote $\sigma(x)=\frac{\omega(x)-|x|}{2}$
\cite{Biane1998,IvanovOlshanski2002} and consider the family of maps
\begin{equation} \label{eq:p-tylda}
\tilde{p}_n(\lambda)= \int_{\R} x^n \sigma''(x) dx.
\end{equation}
Since $\sigma''$ makes sense as a distribution and $\sigma$ is
compactly supported hence the collection
$\big(\tilde{p}_n(\lambda)\big)_n$ determines the Young diagram
$\lambda$ uniquely.



\subsection{Transition measure of a Young diagram}
\index{transition measure of a Young diagram}
\label{subsec:transitionanalytic} To any generalized Young diagram
$\omega$ we can assign the unique probability measure $\mu^{\omega}$
on $\R$, called transition measure of $\omega$, the Cauchy transform
of which
\begin{equation}
\label{eq:Cauchy} G_{\mu^{\omega}}(z)= \int_{\R} \frac{1}{z-x} d
\mu^{\omega}(x)
\end{equation}
is given by
\begin{equation} \label{eq:definicja3} \log G_{\mu^{\omega}}(z)= -\frac{1}{2}
\int_{\R} \log(z-x)  \omega''(x) dx = -\frac{1}{2} \int_{\R}
\frac{1}{z-x} \omega'(x) dx
\end{equation}
%
%
for every $z\notin\R$. For a Young diagram $\lambda$ we will write
$\mu^{\lambda}$ as a short hand of $\mu^{\omega^{\lambda}}$. This
definition may look artificial but it turns out
\cite{Kerov1993transition,OkounkovVershik1996,Biane1998,
Okounkov2000randompermutations} that it is equivalent to natural
representation-theoretic definitions which arise by studying the
representation $\rho_q$ together with the inclusion $S_q\subset
S_{q+1}$.

For $p>0$ and a Young diagram $\lambda$ we consider the rescaled
(generalized) Young diagram $\omega^{p \lambda}$ given by $\omega^{p
\lambda}:x \mapsto p \omega^{\lambda}\big( \frac{x}{p} \big)$.
Informally speaking, the symbol $p \lambda$ corresponds to the shape
of the Young diagram $\lambda$ geometrically scaled by factor $p$
(in particular, if $\lambda$ has $q$ boxes then $p \lambda$ has $p^2
q$ boxes). It is easy to see that \eqref{eq:definicja3} implies that
the corresponding transition measure $\mu^{p \lambda}$ is a dilation
of $\mu^{\lambda}$:
\begin{equation}
\label{eq:skalowanie} \mu^{p \lambda}= D_p \mu^{\lambda}.
\end{equation}
This nice behavior of the transition measure with respect to
rescaling of Young diagrams makes it a perfect tool for the study of
the asymptotics of symmetric groups $S_q$ as $q\to\infty$.

\subsection{Free cumulants of the transition measure}
Cauchy transform of a compactly supported probability measure is
given at the neighborhood of infinity by a power series
$$ G_{\mu}(z)=\frac{1}{z}+\sum_{n\geq 2} M_n z^{-n-1},$$
where $M_n=\int_R x^n d\mu$ are the moments of the measure $\mu$. It
follows that on some neighborhood of infinity $G_{\mu}$ has a right
inverse $K_{\mu}$ with respect to the composition of power series
given by
$$ K_{\mu}(z)=\frac{1}{z}+\sum_{n\geq 1} R_n z^{n-1}$$
convergent on some neighborhood of $0$. The coefficients
$R_i=R_i(\mu)$ are called free cumulants of measure $\mu$. Free
cumulants appeared implicitly in Voiculescu's $R$--transform
\cite{VoiculescuAddition} and their combinatorial meaning was given
by Speicher \cite{Speicher1997}.

Free cumulants are homogenous in the sense that if $X$ is a random
variable and $c$ is some number then
$$ R_i(cX)=c^i R_i(X) $$
and for this reason they are very useful in the study of asymptotic
questions.

Each free cumulant $R_n$ is a polynomial in the moments
$M_1,M_2,\dots,M_n$ of the measure and each moment $M_n$ can be
expressed as a polynomial in the free cumulants $R_1,\dots,R_n$; in
other words the sequence of moments $M_1,M_2,\dots$ and the sequence
of free cumulants $R_1,R_2,\dots$ contain the same information about
the probability measure. The functionals of Young diagrams
considered in \eqref{eq:p-tylda} have a nice geometric
interpretation but they are not very convenient in actual
calculations. For this reason we will prefer to describe the shape
of a Young diagram by considering a family of functionals
\begin{equation}
\label{eq:wolnekumulanty} \lambda\mapsto R_n(\mu^{\lambda})
\end{equation}
given by the free cumulants of the transition measure. Equation
\eqref{eq:definicja3} shows that functionals $\tilde{p}_k$ from the
family \eqref{eq:p-tylda} can be expressed as polynomials in the
functionals from the family \eqref{eq:wolnekumulanty} and vice
versa.

Please note that the first two cumulants of a transition measure do
not carry any interesting information since
$$ R_1(\mu^{\lambda})=M_1(\mu^{\lambda})=0,$$
$$ R_2(\mu^{\lambda})=M_2(\mu^{\lambda})=q,$$
where $q$ denotes the number of the boxes of the Young diagram
$\lambda$.

Above we treated the free cumulant $R_i$ as a function on Young
diagrams, but it also can be viewed (via Fourier transform) as a
central element in $\C(S_q)$.

\section{Representations with character factorization property}
\label{sec:representations}

\subsection{Factorization of characters and Gaussian fluctuations}

The following theorem is the main result of this article. In order
not to scare the Reader we postpone its proof to Section
\ref{sec:Proof-main-result}.

\begin{theoremanddefinition}
\label{theo:main} For each $q\geq 1$ let $\rho_q$ be a
representation of $S_q$. We say that the sequence $(\rho_q)$
has the character factorization property if it fulfills one (hence
all) of the following equivalent conditions:
\begin{itemize}

\item for any cycles $\sigma_1,\dots,\sigma_n$ with disjoint
supports
\begin{equation}
\label{eq:faktoryzacja1} k(\sigma_1,\dots,\sigma_n) q^{
\frac{|\sigma_1|+\cdots+|\sigma_n|+2(n-1)}{2} }= O(1);
\end{equation}

\item for any integers $l_1,\dots,l_n\geq 1$
\begin{equation}
\label{eq:faktoryzacja2}
k^{\bullet}(\Sigma_{l_1},\dots,\Sigma_{l_n})  q^{-
\frac{l_1+\cdots+l_n-n+2}{2}  } = O(1);
\end{equation}

\item for any integers $l_1,\dots,l_n\geq 1$
\begin{equation}
\label{eq:faktoryzacja3} k(\Sigma_{l_1},\dots,\Sigma_{l_n})
q^{-\frac{l_1+\cdots+l_n-n+2}{2} }  = O(1);
\end{equation}

\item for any integers $l_1,\dots,l_n\geq 2$
\begin{equation}
\label{eq:faktoryzacja4} k(R_{l_1},\dots,R_{l_n} )
q^{-\frac{l_1+\cdots+l_n-2(n-1)}{2}  }= O(1) .
\end{equation}
\end{itemize}
\end{theoremanddefinition}

\begin{remark}
In Corollary \ref{cor:roznegradacjejednadegradacja} we will prove
that if conditions \eqref{eq:faktoryzacja2} and
\eqref{eq:faktoryzacja3} hold true then they also hold true in a
more general situation when conjugacy classes $(\Sigma_{l_i})$ with
only one non-trivial cycle are replaced by general conjugacy classes
$\Sigma_{l_{i,1},\dots,l_{i,m(i)}}$. Similarly one can show that if
condition \eqref{eq:faktoryzacja1} holds true then it also holds
true in a general situation when we do not assume that $(\sigma_i)$
are cycles.
\end{remark}

To show that a given sequence of representations has the character
factorization property usually it is the most convenient to verify
condition \eqref{eq:faktoryzacja1} or condition
\eqref{eq:faktoryzacja2}. Then conditions \eqref{eq:faktoryzacja3}
and \eqref{eq:faktoryzacja4} are important corollaries (for
applications see Corollary \ref{coro:gaussian} below).

Expressions appearing in conditions
\eqref{eq:faktoryzacja1}--\eqref{eq:faktoryzacja4} are closely
related to each other and knowledge of one of them allows us to
compute the others (in fact this is how Theorem and Definition
\ref{theo:main} will be proved). For general $n$ these formulas are
quite involved, however in the following we will need only such
formulas for $n\in\{1,2\}$ and these are provided by the following
theorem.

\begin{theorem}
\label{theo:main2} Let $(\rho_q)$ has the character factorization
property. If the limit of one of the expressions
\eqref{eq:faktoryzacja1}--\eqref{eq:faktoryzacja4} exists for
$n\in\{1,2\}$ then the limits of all of the expressions
\eqref{eq:faktoryzacja1}--\eqref{eq:faktoryzacja4} exist for
$n\in\{1,2\}$.

These limits fulfill
\begin{equation}
\label{eq:granicamomentow} c_{l+1}:= \lim_{q\to\infty} \E(\sigma)
q^{ \frac{l-1}{2} }= \lim_{q\to\infty} \E(\Sigma_{l})
q^{-\frac{l+1}{2} }  = \lim_{q\to\infty} \E(R_{l+1})
q^{-\frac{l+1}{2} },
\end{equation}
where $\sigma$ is a cycle of length $l$ and
\begin{multline}
\label{eq:granicamomentow2} \lim_{q\to\infty}
\Cov(R_{l_1+1},R_{l_2+1} ) q^{-\frac{l_1+l_2}{2}}
=\lim_{q\to\infty}\Cov(\Sigma_{l_1},\Sigma_{l_2})
q^{-\frac{l_1+l_2}{2} }=\\
\shoveleft{ \lim_{q\to\infty} \Cov^{\bullet}
(\Sigma_{l_1},\Sigma_{l_2}) q^{-
\frac{l_1+l_2}{2}  } }+ \\ \shoveright{ \sum_{r\geq 1}\sum_{\substack{a_1,\dots,a_r\geq 1\\
a_1+\dots+a_r=l_1}}  \sum_{\substack{b_1,\dots,b_r\geq 1\\
b_1+\dots+b_r=l_2}} \frac{l_1 l_2}{r} c_{a_1+b_1} \cdots
c_{a_r+b_r}=} \\
\shoveleft{ \lim_{q\to\infty} \Cov(\sigma_1,\sigma_2)
q^{\frac{l_1+l_2}{2} } -l_1 l_2 c_{l_1+1} c_{l_2+1}+} \\
\sum_{r\geq 1}\sum_{\substack{a_1,\dots,a_r\geq 1\\
a_1+\dots+a_r=l_1}}  \sum_{\substack{b_1,\dots,b_r\geq 1\\
b_1+\dots+b_r=l_2}} \frac{l_1 l_2}{r} c_{a_1+b_1} \cdots c_{a_r+b_r}
\end{multline}
where $\sigma_1,\sigma_2$ are disjoint cycles of length $l_1,l_2$,
respectively, and where the numbers $c_i$ were defined in
\eqref{eq:granicamomentow}. \end{theorem} Proof of this theorem is
also postponed to Section \ref{sec:Proof-main-result}; we will prove
it together with Theorem \ref{theo:main}. Identity
\eqref{eq:granicamomentow} was proved by Biane
\cite{Biane1998,Biane2001approximate} and we skip its proof; for
Readers acquainted with the results of Section
\ref{sec:Proof-main-result} and Section
\ref{sec:Proof-Theorem-?????} it will be a simple exercise.


\begin{corollary}
\label{coro:gaussian} Let $(\rho_q)$ be as in Theorem
\ref{theo:main2} and let $\lambda$ be a random Young diagram
distributed according to the canonical probability measure
associated to $\rho_q$.
\begin{enumerate}
\item
\label{item:gaba1}\emph{(Gaussian fluctuations of free cumulants)}
Then the joint distribution of the centered random variables
$$r_i=q^{-\frac{i-2}{2}} (R_i - \E R_i)  $$
converges to a Gaussian distribution in the weak topology of
probability measures, where $R_i$ denotes the free cumulant of the
transition measure $\mu^{\lambda}$.
\item \label{item:gaba2} \emph{(Gaussian fluctuations of characters)}
Let $\sigma_i$ denote a cycle of length $i$. Then the joint
distribution of the centered random variables
\begin{equation}
\label{eq:gausowskiecharaktery} q^{ \frac{|\sigma_i|+1}{2} }
\big(\chi^{\lambda}(\sigma_i)- \E \chi^{\lambda}(\sigma_i)\big)
\end{equation}
converges to a Gaussian distribution in the weak topology of
probability measures.
\item \label{item:gaba3}
\emph{(Gaussian fluctuations of the shape of the Young diagrams)}
Then the joint distribution of the centered random variables
$$q^{-\frac{i-2}{2}} (\tilde{p}_i - \E \tilde{p}_i)  $$
converges to a Gaussian distribution in the weak topology of
probability measures, where $\tilde{p}_i=\tilde{p}_i(\lambda)$ is
the functional of the shape of the Young diagram defined in
\eqref{eq:p-tylda}.
\end{enumerate}
\end{corollary}
\begin{proof}
We will prove now point \eqref{item:gaba1}. Condition
\eqref{eq:faktoryzacja4} implies that if $n\neq 2$ then
$$ \lim_{q\to\infty} k(r_{i_1},\dots,r_{i_n})=0 $$
and therefore the family $(r_i)$ converges in moments to a Gaussian
distribution. Since Gaussian measures are uniquely determined by
their moments it follows that the convergence holds true also in the
weak topology of probability measures.

To proof point \eqref{item:gaba2} we observe that asymptotically, as
$q\to\infty$ random variables \eqref{eq:gausowskiecharaktery} have
the same behavior as random variables
\begin{equation}
\label{eq:gausowskiecharaktery2} \lambda\mapsto q^{
\frac{-|\sigma_i|+1}{2} } \big(q^{\underline{|\sigma_i|}}
\chi^{\lambda}(\sigma_i)- \E q^{\underline{|\sigma_i|}}
\chi^{\lambda}(\sigma_i)\big),
\end{equation}
where $q^{\underline{k}}=q (q-1) \cdots (q-k+1)$ denotes the falling
power. Function on Young diagrams \eqref{eq:gausowskiecharaktery2}
corresponds (via Fourier transform) to the central function in
$\C(S_q)$
$$ q^{
\frac{-|\sigma_i|+1}{2} } \big(\Sigma_i - \E \Sigma_i\big).
$$
It follows that we may use \eqref{eq:faktoryzacja3} in the same way
as in the above proof of point \eqref{item:gaba1}.

Ivanov and Olshanski \cite{IvanovOlshanski2002} proved that point
\eqref{item:gaba1} implies point \eqref{item:gaba3}; their proof is
a careful analysis of the fact that $\tilde{p}_i$ can be expressed
as a polynomial in  free cumulants.
\end{proof}

\subsection{Examples}
All examples presented in this section not only have the character
factorization property but additionally are as in Theorem
\ref{theo:main2}.

\begin{example}[Left-regular representation]
\label{example:1} It is easy to check that if $\rho_q$ is the
left-regular representation of $S_q$ then for any permutations
$\sigma_1,\dots,\sigma_n$ with disjoint supports
$$k(\sigma_1,\dots,\sigma_n) q^{
\frac{|\sigma_1|+\cdots+|\sigma_n|+2(n-1)}{2} } = \begin{cases} 1 &
\text{if } n=1 \text{ and } \sigma_1=e, \\ 0 & \text{otherwise.}
\end{cases} $$
It follows from condition \eqref{eq:faktoryzacja1} that the
left-regular representation has the character factorization property
and that the mean and the covariance of the free cumulants are given
by
\begin{equation} \label{eq:kerov1} \lim_{q\to\infty} \E(R_{l+1})
q^{-\frac{l+1}{2} } = \begin{cases} 1 & \text{if } l=1 , \\
0 & \text{if } l\geq 2, \end{cases} \end{equation} and
\begin{equation} \label{eq:kerov2} \lim_{q\to\infty} \Cov(R_{l_1+1},R_{l_2+1} )
q^{-\frac{l_1+l_2}{2}} = \begin{cases} l_1 & \text{if } l_1=l_2\geq
2, \\ 0 & \text{if } l_1=l_2=1, \\ 0 & \text{if } l_1\neq l_2.
\end{cases}
\end{equation}
By applying Corollary \ref{coro:gaussian} we recover Kerov central
limit theorem for the Plancherel measure
\cite{Kerov1993gaussian,IvanovOlshanski2002}.
\end{example}

\begin{example}[Tensor representations]
\label{example:2} For some integer $d_q\geq 1$ let $\rho_q$ be the
representation of $\Sn{q}$ acting on $(\C^{d_q})^{\otimes q}$ by
permutation of factors. This representation appears naturally within
Schur-Weyl duality. It is easy to check that for permutations
$\sigma_1,\dots,\sigma_n$ with disjoint supports
$$k(\sigma_1,\dots,\sigma_n) q^{
\frac{|\sigma_1|+\cdots+|\sigma_n|+2(n-1)}{2} }= \begin{cases}
\left( \frac{\sqrt{q}}{d_q}\right)^{|\sigma_1|} & \text{if } n=1,
\\ 0 & \text{otherwise,} \end{cases} $$
hence if the limit ${p}:=\lim_{q\to\infty} \frac{\sqrt{q}}{d_n}$
exists then (condition \eqref{eq:faktoryzacja1}) the sequence
$(\rho_q)$ has the character factorization property and the mean and
the covariance of the free cumulants are given by
$$ \lim_{q\to\infty} \E(R_{l+1})
q^{-\frac{l+1}{2} }=  {p}^{l-1}, $$ and
$$ \lim_{q\to\infty}
\Cov(R_{l_1+1},R_{l_2+1} ) q^{-\frac{l_1+l_2}{2}}= \sum_{r\geq 2}
\binom{l_1}{r} \binom{l_2}{r} r \ {p}^{l_1+l_2-2 r}
$$
for all integers $l,l_1,l_2\geq 1$. Note that for ${p}=0$ we recover
the fluctuations of the Plancherel measure.

Some asymptotic results for this representation were proved by Biane
\cite{Biane2001approximate}.
\end{example}

\begin{example}[Irreducible representations]
\label{example:3} Let $c>0$ be a constant and let $(\lambda_q)$ be a
sequence of Young diagrams. We assume that $\lambda_q$ has $q$ boxes
and it has at most $c \sqrt{q}$ rows and columns. Suppose that the
shapes of rescaled Young diagrams $q^{-\frac{1}{2} } \lambda_q$
converge to some limit. The convergence of the shapes of Young
diagrams implies convergence of the free cumulants and it follows
(condition \eqref{eq:faktoryzacja4}) that the sequence
$(\rho^{\lambda_q})$ of the corresponding irreducible
representations has the characters factorization property.

In this example the cumulants \eqref{eq:faktoryzacja3} and
\eqref{eq:faktoryzacja4} vanish for $n\geq 2$ since the Young
diagrams are non-random and the corresponding limits for $n=1$ are
determined by the limit of the shape of the Young diagrams.
\end{example}

The above three examples are the building blocks from which one can
construct some more complex representations with the help of the
operations on representations presented below.

\begin{theorem}[Restriction of representations]
Suppose that the sequence of representations $(\rho_q)$ has the
character factorization property. Let a sequence of integers $(r_q)$
be given, such that $r_q\geq q$ and the limit ${p}=\lim_{q\to\infty}
\frac{q}{r_q} $ exists.

Let $\rho'_q$ denote the restriction of the representation
$\rho_{r_q}$ to the subgroup $S_q\subseteq S_{r_q}$. Then the
sequence $(\rho'_q)$ has the factorization property of characters.
The fluctuations of the free cumulants are determined by
\begin{equation}
c_{l+1}':=\lim_{q\to\infty} \E(R'_{l+1}) q^{-\frac{l+1}{2} }=
{p}^{\frac{l-1}{2}} \lim_{q\to\infty} \E(R_{l+1}) q^{-\frac{l+1}{2}
}= {p}^{\frac{l-1}{2}} c_{l+1}, \label{equ:restrictionA}
\end{equation}
\begin{multline} \lim_{q\to\infty} \Cov(R'_{l_1+1},R'_{l_2+1})
q^{-\frac{l_1+l_2}{2}}= \\ {p}^{\frac{l_1+l_2}{2}} \Bigg[
\lim_{q\to\infty} \Cov(R_{l_1+1},R_{l_2+1}) q^{-\frac{l_1+l_2}{2}} -
 l_1 l_2 c_{l_1+1} c_{l_2+1} \left( {p}^{-1}- 1\right) + \\
\sum_{r\geq 1}\sum_{\substack{a_1,\dots,a_r\geq 1\\
a_1+\dots+a_r=l_1}}  \sum_{\substack{b_1,\dots,b_r\geq 1\\
b_1+\dots+b_r=l_2}} \frac{l_1 l_2}{r} c_{a_1+b_1} \cdots c_{a_r+b_r}
\left(  {p}^{-r}-1 \right) \Bigg] \label{mul:restrictionB}
\end{multline}
for all $l,l_1,l_2\geq 1$, where the quantities $R'_i,c'_i$ concern
the representations $(\rho'_q)$ while $R_i,c_i$ concern the
representations $(\rho_q)$.

In particular, for ${p}=0$ we recover the fluctuations of the
Plancherel measure.
\end{theorem}
\begin{proof}
Notice that for any permutations $\sigma_1,\dots,\sigma_n$ the value
of the cumulant $ k(\sigma_1,\dots,\sigma_n)$ is the same for the
representation $\rho'_q$ of $S_q$ and for the representation
$\rho_{r_q}$ of $S_{r_q}$. It follows that
\begin{multline*} \lim_{q\to\infty} k_{\rho'_q}(\sigma_1,\dots,\sigma_n)
q^{\frac{|\sigma_1|+\cdots+|\sigma_n|+2(n-1)}{2}}= \\
{p}^{\frac{|\sigma_1|+\cdots+|\sigma_n|+2(n-1)}{2}}
\lim_{q\to\infty}k_{\rho_{r_q}}(\sigma_1,\dots,\sigma_n)
(r_q)^{\frac{|\sigma_1|+\cdots+|\sigma_n|+2(n-1)}{2}}
\end{multline*}
Since the original representations $(\rho_q)$ fulfill condition
\eqref{eq:faktoryzacja1} hence restricted representations
$(\rho'_q)$ fulfill \eqref{eq:faktoryzacja1} as well. Equations
\eqref{equ:restrictionA} and \eqref{mul:restrictionB} follow as
special cases for $n=1,2$.
\end{proof}
\begin{remark}
Please notice that the above theorem concerns restrictions of the
form $\rho_{r_q}\big\downarrow^{S_{r_q}}_{S_{q}}$ while it is even
more interesting to ask about the asymptotics of the restrictions of
the form $\rho''_q=\rho_{r_q}\big\downarrow^{S_{r_q}}_{S_{q}\times
S_{r_q-q}}$. A typical question is the following one: let
$[\lambda^{(1)}]\times [\lambda^{(2)}]$ be a random irreducible
component of $\rho''_q$; is it true that the joint distribution of
the free cumulants  $\big(R_i(\lambda^{(r)})\big)_{r\in\{1,2\};
i\geq 2}$ converges after appropriate rescaling to a family of
Gaussian variables? The answer for this question is positive and the
Reader may easily prove it using the methods presented in Section
\ref{sec:Proof-main-result}.
\end{remark}

Such representations are no longer representations of symmetric
groups (but products of such groups) and the original definition of
representations with the approximate factorization property does not
fit into this setup. Nevertheless

In the work of Pittel and Romik \cite{PittelRomik04} some
interesting results concerning restrictions of irreducible
representations corresponding to rectangular Young diagrams were
proved.

\begin{theorem}[Outer product of representations]
\label{theo:outerproduct} Suppose that for\/ $i\in\{1,2\}$ the
sequence of representations $(\rho^{(i)}_q)$ has the character
factorization property. Let sequences of positive integers
$r^{(i)}_q$ be given, such that $r^{(1)}_q+r^{(2)}_q=q$ and the
limits ${p}^{(i)}:=\lim_{q\to\infty} \frac{r^{(i)}}{q} $ exist.

Let $\rho'_q=\rho^{(1)}_{r^{(1)}_q} \circ \rho^{(2)}_{r^{(2)}_q}$
denote the outer product of representations. Then the sequence
$(\rho'_q)$ has the factorization property of characters with
\begin{equation}
\label{eq:outerA} c_{l+1}':=\lim_{q\to\infty} \E(R'_{l+1})
q^{-\frac{l+1}{2} }=  \left({p}^{(1)}\right)^{\frac{l+1}{2}}
c^{(1)}_{l+1} + \left({p}^{(2)}\right)^{\frac{l+1}{2}}
c^{(2)}_{l+1},
\end{equation}
and with an explicit (but involved) covariance of free cumulants.
The appropriate disjoint covariance is given by
\begin{multline*} \lim_{q\to\infty}
\Cov_{\rho'_q}^{\bullet}(\Sigma_{l_1},\Sigma_{l_2})
q^{-\frac{l_1+l_2}{2}}= \\\shoveleft{({p}^{(1)}
)^{\frac{l_1+l_2}{2}} \lim_{q\to\infty}
\Cov_{\rho^{(1)}_{q}}^{\bullet}(\Sigma_{l_1},\Sigma_{l_2})
q^{-\frac{l_1+l_2}{2}}+} \\
({p}^{(2)} )^{\frac{l_1+l_2}{2}}  \lim_{q\to\infty}
\Cov_{\rho^{(2)}_{q}}^{\bullet}(\Sigma_{l_1},\Sigma_{l_2})
q^{-\frac{l_1+l_2}{2}}.
\end{multline*}
\end{theorem}
\begin{proof}
By definition,
$$ \rho'_q=\big( \rho^{(1)}_{r^{(1)}_q} \times
\rho^{(2)}_{r^{(2)}_q} \big) \big\uparrow_{S_{r^{(1)}_q} \times
S_{r^{(2)}_q}}^{S_q}$$ and we can use the Frobenius reciprocity
between the induction and restriction of representations. It follows
that the corresponding normalized characters fulfill for all
$l_1,\dots,l_n\geq 1$
$$ \chi_{\rho'} (\Sigma_{l_1} \bullet \cdots \bullet \Sigma_{l_n})=
\big(\chi_{\rho^{(1)}} \otimes \chi_{\rho^{(2)}} \big) \big(
(\Sigma_{l_1}\otimes 1+1 \otimes \Sigma_{l_1}) \bullet \cdots
\bullet (\Sigma_{l_n}\otimes 1+1 \otimes \Sigma_{l_n})  \big).$$ We
can treat the left-hand side as a mixed moment $\E (\Sigma_{l_1}
\bullet \cdots \bullet \Sigma_{l_n})$ in the algebra of conjugacy
classes equipped with the disjoint product; analogous interpretation
is possible also for the right-hand side. Since equality holds for
all $l_1,\dots,l_n\geq 2$ it follows that the joint distributions of
the family of random variables $(\Sigma_i)_{i\geq 2}$ coincides with
the joint distribution of random variables $(\Sigma_{i}\otimes 1+1
\otimes \Sigma_{i})_{i\geq 2}$; in particular their cumulants are
equal:
\begin{multline}
\label{eq:kumulantyok}
k^{\bullet}_{\rho'}(\Sigma_{l_1},\dots,\Sigma_{l_n})=
k^{\bullet}_{\rho^{(1)}\otimes \rho^{(2)}}(\Sigma_{l_1}\otimes
1+1\otimes \Sigma_{l_1},\dots,\Sigma_{l_n}\otimes 1+1\otimes
\Sigma_{l_n})=\\
k^{\bullet}_{\rho^{(1)}}(\Sigma_{l_1},\dots,\Sigma_{l_n})+
k^{\bullet}_{\rho^{(2)}}(\Sigma_{l_1},\dots,\Sigma_{l_n}).
\end{multline}
In the last equality we used that the cumulant is linear with
respect to each of the arguments and that the mixed cumulant of
independent random variables vanishes.

It follows that \begin{multline*} \lim_{q\to\infty}
k_{\rho'_q}^{\bullet}(\Sigma_{l_1},\dots,\Sigma_{l_n})
q^{-\frac{l_1+\cdots+l_n-n+2}{2}}= \\({p}^{(1)}
)^{\frac{l_1+\cdots+l_n-n+2}{2}}  \lim_{q\to\infty}
k_{\rho^{(1)}_{r^{(1)}_q}}^{\bullet}(\Sigma_{l_1},\dots,\Sigma_{l_n})
(r^{(1)}_q)^{-\frac{l_1+\cdots+l_n-n+2}{2}}+ \\
({p}^{(2)} )^{\frac{l_1+\cdots+l_n-n+2}{2}}  \lim_{q\to\infty}
k_{\rho^{(2)}_{r^{(2)}_q}}^{\bullet}(\Sigma_{l_1},\dots,\Sigma_{l_n})
(r^{(1)}_q)^{-\frac{l_1+\cdots+l_n-n+2}{2}}.
\end{multline*}
It follows that if the representations $\rho^{(i)}$ fulfill
condition \eqref{eq:faktoryzacja2} then representations $\rho'_q$
fulfill \eqref{eq:faktoryzacja2} as well. Equation \eqref{eq:outerA}
results as a special case of $n=1$ and we leave it as a simple
exercise to the Reader to study the case $n=2$ and to find the
covariance of the free cumulants, analogous to
\eqref{equ:restrictionA}.
\end{proof}

\begin{theorem}[Induction of representations]
Suppose that the sequence of representations $(\rho_q)$ has
character factorization property. Let a sequence of integers $r_q$
be given, such that $r_q\leq q$ and the limit ${p}=\lim_{q\to\infty}
\frac{r_q}{q} $ exists.

Let $\rho'_q=\rho_{r_q}\uparrow_{S_{r_q}}^{S_q}$ denote the induced
representation. Then the sequence $(\rho'_q)$ has the characters
factorization property  with
$$ c'_{l+1}=\begin{cases} {p}^{\frac{l+1}{2}} c_{l+1} & \text{for } l\geq 2, \\
1 & \text{for } l=1, \end{cases} $$ and with an explicit (but
involved) covariance of free cumulants.
\end{theorem}
\begin{proof}
It is enough to adapt the proof of Theorem \ref{theo:outerproduct}.
\end{proof}

\begin{theorem}[Tensor product of representations]
Suppose that for $i\in\{1,2\}$ the sequence of representations
$(\rho^{(i)}_q)$ has character factorization property. Then the
tensor product $\rho'_q=\rho^{(1)}_q \otimes \rho^{(2)}_q$ has the
property of factorization of characters. Furthermore, the limit
distribution and the fluctuations are the same as for the Plancherel
measure \eqref{eq:kerov1} and \eqref{eq:kerov2}.
\end{theorem}
\begin{proof}
Since the normalized characters fulfill for any $\pi\in S_q$
$$ \chi_{\rho'_q}(\pi)= \big( \chi_{\rho^{(1)}_q} \otimes
\chi_{\rho^{(2)}_q} \big) (\pi \otimes \pi)$$ hence also the
corresponding cumulants are equal:
$$k_{\rho'_q}(\sigma_1,\dots,\sigma_n)=k_{\rho^{(1)}_q\otimes \rho^{(2)}_q}
(\sigma_1\otimes \sigma_1,\dots,\sigma_n\otimes \sigma_n).$$ Theorem
\ref{theo:leonovsiraev} of Leonov and Sirjaev  can be used to
calculate the right-hand side. Lemma \ref{lem:connections} together
with the condition \eqref{eq:faktoryzacja1} for $\rho^{(i)}$ show
that the right-hand side is of order
$O(q^{-(|\sigma_1|+\cdots+|\sigma_n|+n-1)})  $ hence condition
\eqref{eq:faktoryzacja1} is fulfilled for $\rho'$. It also follows
that the limits in Theorem \ref{theo:main2} are given by
$$c_{l+1}=\lim_{q\to\infty} \E(\sigma)
q^{ \frac{l-1}{2} }= \begin{cases} 1 & \text{if } l=1, \\ 0 &
\text{if } l\geq 2, \end{cases}
$$
$$ \lim_{q\to\infty} \Cov(\sigma_1,\sigma_2)
q^{\frac{l_1+l_2}{2} }= 0, $$ where $\sigma$, $\sigma_1$, $\sigma_2$
are disjoint cycles with lengths, respectively, $l$, $l_1$, $l_2$.
\end{proof}

\section{Proof of the main result}
\label{sec:Proof-main-result}

\subsection{Toy example}

Let us have a look on the main result (Theorem and Definition
\ref{theo:main}) for the simplest nontrivial case of $n=2$. This
will give us a heuristical insight into the problems which we shall
encounter in the proof.

As we shall see in the following, the proof of the equivalence
\eqref{eq:faktoryzacja3} and \eqref{eq:faktoryzacja4} is very easy,
therefore the main difficulty is to show the equivalence of the
conditions \eqref{eq:faktoryzacja1}, \eqref{eq:faktoryzacja2},
\eqref{eq:faktoryzacja3}. For simplicity, we shall concentrate on
the implications in only one direction
\eqref{eq:faktoryzacja1}$\implies$\eqref{eq:faktoryzacja2}$\implies$\eqref{eq:faktoryzacja3}.

The condition \eqref{eq:faktoryzacja3} for $n=2$ requires that
$k(\Sigma_{l_1},\Sigma_{l_2})$ should not grow too fast. The
identity
\begin{multline}
\label{eq:problemy} k(\Sigma_{l_1},\Sigma_{l_2})=\chi(\Sigma_{l_1}
\Sigma_{l_2})-\chi(\Sigma_{l_1}) \chi(\Sigma_{l_2}) =\\
\chi(\Sigma_{l_1} \Sigma_{l_2}-\Sigma_{l_1}\bullet \Sigma_{l_2})+
\big[ \chi(\Sigma_{l_1} \bullet \Sigma_{l_2}) - \chi(\Sigma_{l_1})
\chi(\Sigma_{l_2}) \big]
\end{multline}
shows that there are three reasons why
$k(\Sigma_{l_1},\Sigma_{l_2})$  is non-zero:
\begin{enumerate}
\item \label{diff:1} the difference $\Sigma_{l_1} \Sigma_{l_2}-\Sigma_{l_1}\bullet
\Sigma_{l_2}$ is non-zero. We need to estimate the conjugacy classes
contributing to this difference.
\item \label{diff:2} there are $q\podpotega{l_1} q\podpotega{l_2}$ summands which
contribute to $\chi(\Sigma_{l_1}) \chi(\Sigma_{l_2})$ while there
are only $q\podpotega{l_1+l_2}$ summands which contribute to
$\chi(\Sigma_{l_1} \bullet\Sigma_{l_2})$. Under the simplifying
assumption that $\chi(\pi_1 \pi_2)\approx \chi(\pi_1) \chi(\pi_2)$
the second summand in \eqref{eq:problemy} is therefore of order
$$(q\podpotega{l_1+l_2}-q\podpotega{l_1} q\podpotega{l_2}) \chi(\pi_1
\pi_2).$$ We need to find an estimate for
$(q\podpotega{l_1+l_2}-q\podpotega{l_1} q\podpotega{l_2})$.
\item \label{diff:3} every summand contributing to $\chi(\Sigma_{l_1})
\chi(\Sigma_{l_2})$ is equal to $\chi(\pi_1 \pi_2)$ while every
summand contributing to $\chi(\Sigma_{l_1,l_2})$ is equal to
$\chi(\pi_1) \chi(\pi_2)$. We need to find an estimate for
$(\chi(\pi_1 \pi_2)-\chi(\pi_1) \chi(\pi_2))$.
\end{enumerate}
The difficulty caused by \eqref{diff:3} can be very easily overcome:
it is basically the condition \eqref{eq:faktoryzacja1}. Our proof of
the main result will be therefore divided into two parts; each
devoted to one of the remaining difficulties.


Note that the second summand in \eqref{eq:problemy} can be written
as
$$k^{\bullet}(\Sigma_{l_1},\Sigma_{l_2})=\chi(\Sigma_{l_1} \bullet \Sigma_{l_2}) - \chi(\Sigma_{l_1})
\chi(\Sigma_{l_2});$$ in other words the proof of the proof of the
implication
\eqref{eq:faktoryzacja1}$\implies$\eqref{eq:faktoryzacja2} is
equivalent to overcoming the difficulty \eqref{diff:2} and the proof
of the implication
\eqref{eq:faktoryzacja2}$\implies$\eqref{eq:faktoryzacja3} is
equivalent to overcoming the difficulty \eqref{diff:1}.

\subsection{Partitions} \label{subsec:partitions} \index{partition}
We recall that $\pi=\{\pi_1,\dots,\pi_r\}$ is a partition of a
finite ordered set $X$ if sets $\pi_1,\dots,\pi_r$ are nonempty and
disjoint and if $\pi_1\cup\cdots\cup\pi_r=X$.


We say that a partition $\pi$ is smaller (or finer) than a partition
$\rho$ of the same set if every block of $\pi$ is a subset of some
block of $\rho$ and we denote it by \index{$Pi$@$\pi\leq\rho$}
$\pi\leq\rho$. For partitions $\pi_1,\pi_2$ of the same set we
denote by $\pi_1\vee \pi_2$ the minimal partition which is greater
or equal than both $\pi_1$ and $\pi_2$. The set of all partitions of
some given set $X$ has a maximal element, namely the partition with
only one block equal to $X$.

\subsection{Algebra of conjugacy classes and its filtration}
\label{subsec:algebraofconjugacy} Ivanov and Kerov
\cite{IvanovKerov1999} defined a partial permutation of a set $X$ as
a pair $\alpha=(d,w)$, where $d$ (called support of $\alpha$) is any
subset of $X$ and $w:X\rightarrow X$ is a bijection which equal to
identity outside of $d$. The usual product of partial permutations
is given by
$$ (d_1,w_1)(d_2,w_2)=(d_1\cup d_2,w_1 w_2)$$
and their disjoint product $\bullet$ was defined in
\eqref{eq:definition-of-disjoint}. Partial permutations behave like
the usual permutations (to a partial permutation we can canonically
associate the usual permutation $w$) except that we can distinguish
two kinds of fix-points $x$ for a partial permutations: true
fix-points (i.e.\ $x\notin d$) and cycles of length one ($x\in d$,
$w(x)=x$). Partial permutations form a semigroup; in this article we
are interested also in the corresponding semigroup algebra which
should be regarded as an analogue of the permutation group algebra
$\C(S_q)$ equipped with some additional structure.

In fact, to define correctly the notion of the disjoint product
$\bullet$ from Section \ref{subsec:disjoint} we must use the
semigroup algebra corresponding to partial permutations and not the
group algebra $\C(S_q)$. The reason for this is that we must
distinguish two kinds of fix-points since, for example,
$\Sigma_{1,1}$ and $\Sigma_1$ represent multiples of each other in
$\C(S_q)$, but the disjoint product treats them differently.

%
%

One can show
\cite{IvanovKerov1999,SniadyGaussianFluctuationsAndGenus} that the
family of normalized conjugacy classes $(\Sigma_{k_1,\dots,k_n})$
and the family of free cumulants $(R_i)$ generate the same filtered
algebra, called algebra of conjugacy classes, when for the degrees
of the generators we take
$$ \deg \Sigma_{k_1,\dots,k_n}=(k_1+1)+\cdots+(k_n+1),$$
$$ \deg R_i=i.  $$
The above statement holds true both when as the product we take the
usual product of partial permutations (in this case we denote the
resulting algebra by $\A$) and when as the product we take the
disjoint product (in this case we denote the resulting algebra by
$\A^{\bullet}$).

One can also show
\cite{IvanovKerov1999,SniadyGaussianFluctuationsAndGenus} that in
the scaling considered in this article (i.e.\ a typical Young
diagram has at most $O(\sqrt{q})$ rows and columns) so defined
degree determines the asymptotic behavior of an element, namely
$$ \E X = O\big( q^{\frac{\deg X}{2}}\big) $$
for any element $X$ of this algebra. Furthermore, the generators
fulfill
\begin{equation}
\label{eq:wolnekumulantyiklasy} R_i= \Sigma_{i-1} + \text{(terms of
degree at most $i-2$)}.
\end{equation}

\subsection{Three probability spaces}
Commutative algebras $\A$ and $\A^{\bullet}$ can be regarded as
algebras of random variables on a purely formal level (usually it is
not possible to represent them as algebras of functions on some
Kolmogorov probability space).

Algebras $\A$ and $\A^{\bullet}$ are trivially isomorphic as vector
spaces; we will denote by $\E^{\text{id}}:\A\rightarrow\A^{\bullet}$
the identity map between them, in other words
$$ \E^{\text{id}}(x)=x. $$
One can think that $\E^{\text{id}}$ is a kind of a conditional
expectation.

If $\rho$ is a representation of $S_q$ we consider maps
$\E:\A\rightarrow\C$ and $\E^{\bullet}:\A^{\bullet}\rightarrow\C$
given by
$$ \E(x)=\E^{\bullet}(x)=\chi_{\rho}(x). $$

In this way the following diagram commutes:
\begin{equation}
\label{eq:diagram} \xymatrix{ {\A} \ar[r]_-{\E^{\text{id}} }
\ar@/^1pc/[rr]^-{\E} & {\A^{\bullet}} \ar[r]_-{\E^{\bullet}} & {\C}
}.
\end{equation}
and we may consider three different probability structures:
\begin{itemize}
 \item algebra $\A$
equipped with the expectation $\E$ (which gives rise to the natural
cumulants $k$),
 \item algebra $\A^{\bullet}$ equipped with the
expectation $\E^{\bullet}$ (which gives rise to the disjoint
cumulants $k^{\bullet}$),
 \item algebra $\A$ equipped with the conditional
expectation $\E^{\text{id}}:\A\rightarrow\A^{\bullet}$
(corresponding conditional cumulants belong to $\A^{\bullet}$ and
will be denoted by $k^{\text{id}}$).
\end{itemize}

The commutativity of the diagram \eqref{eq:diagram} implies that the
relation between the corresponding three cumulants is given by the
following formula of Brillinger \cite{Brillinger69} (see also
\cite{Lehner04}).
\begin{proposition}
\label{prop:Brillinger} If\/ $x_1,\dots,x_n$ belong to the center
of\/ $\A$ then \begin{equation} k(x_1,\dots,x_n)=\sum_{\pi}
k^{\bullet}\big[ k^{\text{\rm{id}}}(x_i: i\in \pi_j   ): j=1,2,\dots
\big] , \label{eq:Brillinger} \end{equation}
 where the sum runs over all partitions
$\pi=\{\pi_1,\pi_2,\dots\}$ of the set $\{1,\dots,n\}$.
\end{proposition}

\begin{example}
Let us consider $n=2$: there are two partitions of $\{1,2\}$, namely
$\big\{ \{1\}, \{2\} \big\}$ and $\big\{ \{1,2\} \big\}$ therefore
$$ k(x_1,x_2)= k^{\bullet}\big( k^{\text{\rm{id}}}(x_1), k^{\text{\rm{id}}}(x_2) \big)+
k^{\bullet}\big( k^{\text{\rm{id}}}(x_1,x_2) \big).
  $$
Similarly, for $n=3$ \begin{multline*} k(x_1,x_2,x_3)=
k^{\bullet}\big( k^{\text{\rm{id}}}(x_1),
k^{\text{\rm{id}}}(x_2),k^{\text{\rm{id}}}(x_2) \big)+
k^{\bullet}\big( k^{\text{\rm{id}}}(x_1,x_2),k^{\text{\rm{id}}}(x_3)
\big)+\\ k^{\bullet}\big(
k^{\text{\rm{id}}}(x_1,x_3),k^{\text{\rm{id}}}(x_2) \big)+
k^{\bullet}\big( k^{\text{\rm{id}}}(x_2,x_3),k^{\text{\rm{id}}}(x_1)
\big)+ \\ k^{\bullet}\big( k^{\text{\rm{id}}}(x_1)
,k^{\text{\rm{id}}}(x_2),k^{\text{\rm{id}}}(x_3) \big).
\end{multline*}

\end{example}

The following result will be of great importance in this article.
\begin{theorem}
\label{theo:kumulantykombinatoryczne} For any $x_1,\dots,x_n\in\A$
\begin{equation}
\label{eq:kumulantykombinatoryczne} \deg k^{\text{\rm
id}}(x_1,\dots,x_n) \leq \deg x_1+\cdots+ \deg x_n-2(n-1).
\end{equation} Furthermore, the highest-order term of the second
cumulant is given by
\begin{multline}
\label{eq:a-tak-dokladnie?} k^{\text{\rm
id}}(\Sigma_{l_1},\Sigma_{l_2})=
\sum_{r\geq 1}\sum_{\substack{a_1,\dots,a_r\geq 1\\
a_1+\dots+a_r=l_1}}  \sum_{\substack{b_1,\dots,b_r\geq 1\\
b_1+\dots+b_r=l_2}} \frac{l_1 l_2}{r}
\Sigma_{(a_1+b_1-1),\dots,(a_r+b_r-1)}
+ \\
\text{terms of degree at most $(l_1+l_2-2)$}.
\end{multline}
\end{theorem}
We postpone its proof to Section \ref{sec:Proof-Theorem-?????}.

\subsection{Cumulants of products}
If $\rho=\{\rho_1,\dots,\rho_k\}$ is a partition of the set
$\{1,\dots,n\}$ with blocks $\rho_i=\{\rho_{i,1},\dots,\rho_{i,m(i)}
\}$ we define partition--indexed cumulants given by a multiplicative
extension of the usual cumulants:
\begin{equation}
\label{eq:multiplikatywnie} k_{\rho}(X_1,\dots,X_n)= \prod_i
k(X_{\rho_{i,1} }, X_{\rho_{i,2}},\dots, X_{\rho_{i,m(i)}});
\end{equation}
for example
$$  k_{ \{1,3,4\},\{2,5\}} (X_1,X_2,X_3,X_4,X_5)= k (X_1,X_3,X_4)\
k (X_2,X_5). $$

The following formula for cumulants of products of random variables
was proved by Leonov and Sirjaev \cite{LeonovSiraev}.

\begin{theorem}
\label{theo:leonovsiraev} Let $i_1<i_2<\cdots<i_{n+1}$ be integers
and let $X_{i_1+1},X_{i_1+2},\dots,X_{i_{n+1}}$ be a family of
random variables; then
\begin{equation} \label{eq:leonovsiraev} k\Big(\prod_{i_1+1\leq
j\leq i_2} X_{j}, \dots , \prod_{i_n+1\leq j\leq i_{n+1}} X_{j}
\Big) = \sum_{\pi} k_{\pi}(X_{i_1+1},X_{i_1+2},\dots,X_{i_{n+1}}),
\end{equation}
where the sum runs over all partitions $\pi$ of the set
$\{i_1+1,i_1+2,\dots,i_{n+1}\}$ with the additional property that
$$\pi \vee \big\{ \{i_1+1,i_1+2,\dots,i_2\}, \dots,
\{i_{n}+1,i_{n}+2,\dots,i_{n+1} \}   \big\}$$ is the maximal
partition with only one block.
\end{theorem}

\begin{corollary}
\label{cor:takieisiakie} Let permutations $\sigma_1,\dots,\sigma_n$
be disjoint cycles of length $l_1,\dots,l_n$. Then
\begin{multline}
\label{eq:happynewyear}
k^{\bullet}(\Sigma_{l_1},\dots,\Sigma_{l_n})= \\ \sum_{\pi^{(1)},
\pi^{(2)}} k_{\pi^{(1)}}(\sigma_{1},\dots,\sigma_{n}) \
k^{\bullet}_{\pi^{(2)}} (\Sigma_{\underbrace{1,\dots,1}_{l_1 \text{
times}} },\dots,\Sigma_{\underbrace{1,\dots,1}_{l_n \text{ times}}
}),
\end{multline}
where the sum runs over all partitions $\pi^{(1)},\pi^{(2)}$ of the
set $\{1,\dots,n\}$ such that $\pi^{(1)}\vee \pi^{(2)}= \big\{
\{1,\dots,n\} \big\} $.
\end{corollary}

\begin{proof}
We consider the algebra $\C(S_q)\otimes \A^{\bullet}$ with the usual
product on the first factor and the disjoint product on the second
one: $(a_1\otimes b_1) (a_2\otimes b_2)=(a_1 a_2) \otimes (b_1
\bullet b_2)$.

For all integers $l_1,\dots,l_n\geq 1$ and disjoint cycles
$\sigma_1,\dots,\sigma_{n}$ with appropriate lengths we clearly have
$$ \E \big(\Sigma_{l_1} \bullet \cdots \bullet \Sigma_{l_n} \big)=
\E \big[ (\sigma_{1} \otimes \Sigma_{\underbrace{1,\dots,1}_{l_1
\text{ times}} } ) \cdots (\sigma_{n} \otimes
\Sigma_{\underbrace{1,\dots,1}_{l_n \text{ times}} } ) \big],
$$
where the product on the left-hand side is taken in $\A^{\bullet}$
and the product on the right-hand side is taken in $\C(S_q)\otimes
\A^{\bullet}$. In other words: the mixed moments of the family of
random variables $(\Sigma_{l_i})$ coincide with the mixed moments of
the family $(\sigma_{i} \otimes \Sigma_{\underbrace{1,\dots,1}_{l_i
\text{ times}} }) $. It follows that corresponding cumulants are
equal:
$$ k^{\bullet} \big(\Sigma_{l_1},\dots,\Sigma_{l_n} \big)=
k(\sigma_{1} \otimes \Sigma_{\underbrace{1,\dots,1}_{l_1 \text{
times}} } , \dots ,\sigma_{n} \otimes
\Sigma_{\underbrace{1,\dots,1}_{l_n \text{ times}} } ).$$ We use now
Theorem \ref{theo:leonovsiraev} to compute the right-hand side which
finishes the proof.
\end{proof}

In applications of this result we will find useful the following
lemma.
\begin{lemma}
\label{lem:connections} Let $\pi=\{\pi_1,\dots,\pi_m\}$ be a
partition which contributes to the right-hand side of
\eqref{eq:leonovsiraev}. Then
$$ \sum_i \big( |\pi_i|-1 \big) \geq n-1. $$

Let $\pi^{(k)}=\{\pi^{(k)}_1,\dots,\pi^{(k)}_{m^{(k)}}\}$ for
$k\in\{1,2\}$ be partitions which contribute to the right-hand side
of \eqref{eq:happynewyear}. Then
$$ \sum_i \big( |\pi^{(1)}_i|-1 \big)+\sum_i \big( |\pi^{(2)}_i|-1 \big) \geq n-1. $$
\end{lemma}
\begin{proof}
Partition $\pi$ can be obtained from the trivial partition (every
block consists of exactly one element) by performing $\sum_i \big(
|\pi_i|-1 \big)$ times the following operation: we select two blocks
and merge them into a single one. Clearly, the number of blocks of
partition
$$ \pi \vee \big\{ \{i_1+1,i_1+2,\dots,i_2\}, \dots,
\{i_{n}+1,i_{n}+2,\dots,i_{n+1} \}   \big\}$$ either decreases by
one or stays the same in each step. The initial number of blocks is
equal to $n$ and the final number of blocks is equal to $1$, which
finishes the proof.

The proof of the second statement is analogous and we skip it.
\end{proof}


\begin{corollary}
\label{cor:roznegradacjejednadegradacja} Let $X$ be a set of
generators of the algebra of conjugacy classes $\A$. Suppose that
\begin{equation}
\label{eq:quickdecay} k(a_1,\dots,a_n) q^{-\frac{\deg
a_1+\cdots+\deg a_n-2(n-1)}{2}  }= O(1)
\end{equation}
holds true for all $n\geq 1$ and $a_1,\dots,a_n\in X$. Then
\eqref{eq:quickdecay} holds true for all $n\geq 1$ and
$a_1,\dots,a_n\in \A$.
\end{corollary}
\begin{proof}
Clearly, it is enough to consider the case when $a_1,\dots,a_n$ are
monomials in elements of $X$. In order to estimate each summand on
the right-hand side of \eqref{eq:leonovsiraev} we apply Lemma
\ref{lem:connections}.
\end{proof}

\begin{remark}
Note that the above Corollary holds true also when the cumulants $k$
are replaced by $k^{\bullet}$. In particular, this Corollary can be
applied for \eqref{eq:faktoryzacja2} and \eqref{eq:faktoryzacja3}.
\end{remark}


\subsection{Proof of the main theorem: equivalence
   \eqref{eq:faktoryzacja2}$\iff$\eqref{eq:faktoryzacja3}}
\begin{proof}[Proof of the implication
   \eqref{eq:faktoryzacja2} $\implies$ \eqref{eq:faktoryzacja3}]
Our goal is to use Proposition \ref{prop:Brillinger} in order to
express $k^{\bullet}(\Sigma_{l_1},\dots,\Sigma_{l_n})$ in terms of
the cumulants $k^{\text{id}}$ and $k^{\bullet}$. To estimate a
summand on the right-hand side of \eqref{eq:Brillinger}
corresponding to a partition $\pi=\{\pi_1,\dots,\pi_m\}$ of
$\{1,\dots,n\}$ we use Theorem \ref{theo:kumulantykombinatoryczne}
and get
$$\deg k^{\text{\rm{id}}}(x_i: i\in \pi_j)\leq \Big( \sum_{i\in\pi_j} \deg x_i
\Big)-2(|\pi_j|-1).$$ Thus
\begin{equation}
\label{equ:szacowanie} \Big(\sum_{1\leq j\leq m} \deg
k^{\text{\rm{id}}}(x_i: i\in \pi_j)\Big) - 2(m-1) \leq \Big( \sum_i
\deg x_i \Big)-2(n-1).
\end{equation}

Assumption \eqref{eq:faktoryzacja2} and Corollary
\ref{cor:roznegradacjejednadegradacja} show that
$$ k^{\bullet}\big[ k^{\text{\rm{id}}}(x_i: i\in \pi_j   ):
j=1,2,\dots \big] q^{\frac{-\sum_j \deg k^{\text{\rm{id}}}(x_i: i\in
\pi_j) + 2(m-1)}{2}}= O(1).$$ Now \eqref{equ:szacowanie} finishes
the proof.
\end{proof}

As a byproduct, in a similar way we obtain a proof of the identity
\begin{multline}
\label{eq:tttt}
 \lim_{q\to\infty}\Cov(\Sigma_{l_1},\Sigma_{l_2})
q^{-\frac{l_1+l_2}{2} }= \lim_{q\to\infty} \Cov^{\bullet}
(\Sigma_{l_1},\Sigma_{l_2}) q^{-
\frac{l_1+l_2}{2}  } +   \\ \sum_{r\geq 1}\sum_{\substack{a_1,\dots,a_r\geq 1\\
a_1+\dots+a_r=l_1}}  \sum_{\substack{b_1,\dots,b_r\geq 1\\
b_1+\dots+b_r=l_2}} \frac{l_1 l_2}{r} c_{a_1+b_1} \cdots c_{a_r+b_r}
\end{multline}
which is a part of Theorem \ref{theo:main2}.

\begin{proof}[Proof of the implication
   \eqref{eq:faktoryzacja3} $\implies$ \eqref{eq:faktoryzacja2}]
We will use induction with respect to $(\deg
\Sigma_{l_1}+\cdots+\deg \Sigma_{l_n})$. Equation
\eqref{eq:Brillinger} can be written in the form
\begin{multline*}
k^{\bullet}(\Sigma_{l_1},\dots,\Sigma_{l_n}) =
k(\Sigma_{l_1},\dots,\Sigma_{l_n})- \\ \sum_{\pi\neq\{ \{1,\dots,n\}
\}} k^{\bullet}\big[ k^{\text{\rm{id}}}(x_i: i\in \pi_j   ):
j=1,2,\dots \big] .
\end{multline*}
The inductive hypothesis can be used to estimate the right-hand side
in the same way as in the proof of the opposite implication above.
\end{proof}

\subsection{Proof of the main theorem: equivalence \eqref{eq:faktoryzacja1} $\iff$
\eqref{eq:faktoryzacja2}}

\subsubsection{Cumulants of falling factorials}

Element $\Sigma_{\underbrace{1,\dots,1}_{k \text{ times}} }\in
\C(S_q)$ is equal to $q (q-1) \cdots (q+1-k)$, the multiple of
identity, therefore no matter which representation we consider we
always have
$$ \E \Sigma_{\underbrace{1,\dots,1}_{k \text{ times}} }\in
\C(S_q)=q (q-1) \cdots (q+1-k).$$ However, it should be stressed
that if we consider the algebra of conjugacy classes equipped with
the disjoint product then this element is not longer a multiple of
identity and
$$\Sigma_{\underbrace{1,\dots,1}_{k_1 \text{
times}} }\bullet \Sigma_{\underbrace{1,\dots,1}_{k_2 \text{ times}}
} =\Sigma_{\underbrace{1,\dots,1}_{k_1+k_2 \text{ times}} }.$$

\begin{lemma}
\label{lem:kumulantyfallingfactorials} For any integers
$l_1,\dots,l_n\geq 1$
$$ k^{\bullet}(\Sigma_{\underbrace{1,\dots,1}_{l_1 \text{ times}} },\dots,
\Sigma_{\underbrace{1,\dots,1}_{l_n \text{ times}} })=
O(q^{l_1+\cdots+l_n+1-n}). $$
\end{lemma}
\begin{proof}
It is easy to check that
$$k(\Sigma_{\underbrace{1,\dots,1}_{l_1 \text{ times}}
},\dots,\Sigma_{\underbrace{1,\dots,1}_{l_n \text{ times}}
})=\begin{cases} q (q-1) \cdots (q+1-l_1) & \text{if } n=1, \\
0 & \text{if } n\geq 2, \end{cases} $$ hence in particular
$$k(\Sigma_{\underbrace{1,\dots,1}_{l_1 \text{ times}}
},\dots,\Sigma_{\underbrace{1,\dots,1}_{l_n \text{ times}} })
q^{-(l_1+\cdots+l_n+1-n)} = O(1).$$  We leave it as a simple
exercise to the reader to check that the presented above proof of
the implication
\eqref{eq:faktoryzacja3}$\implies$\eqref{eq:faktoryzacja2} can be
applied here.
%
%
%
%
%
%
\end{proof}

\begin{proof}[Proof of the implication
     \eqref{eq:faktoryzacja1} $\implies$ \eqref{eq:faktoryzacja2}]
Our goal is to estimate the right-hand side of
\eqref{eq:happynewyear}. The assumption \eqref{eq:faktoryzacja1}
implies that
$$ k_{\pi_{(1)}}(\sigma_1,\dots,\sigma_n) = O\big(q^{
-\frac{|\sigma_1|+\cdots+|\sigma_n|+2 n-2 \text{(number of blocks of
$\pi_{(1)}$)}}{2} }\big);
$$
Lemma \ref{lem:kumulantyfallingfactorials} implies that
$$ k^{\bullet}_{\pi_{(2)}} (\Sigma_{\underbrace{1,\dots,1}_{l_1 \text{
times}} },\dots,\Sigma_{\underbrace{1,\dots,1}_{l_n \text{ times}}
}) = O\big( q^{l_1+\cdots+l_n +\text{(number of blocks of
$\pi_{(2)}$)}  -n} \big). $$ Now it is enough to use Lemma
\ref{lem:connections}.
\end{proof}

We leave the proof of the identity
$$
 \lim_{q\to\infty} \Cov^{\bullet} (\Sigma_{l_1},\Sigma_{l_2}) q^{-
\frac{l_1+l_2}{2}  } =  \lim_{q\to\infty} \Cov(\sigma_1,\sigma_2)
q^{\frac{l_1+l_2}{2} } -l_1 l_2 c_{l_1+1} c_{l_2+1}
$$
which is a part of Theorem \ref{theo:main2} as a simple exercise.

\begin{proof}[Proof of \eqref{eq:faktoryzacja2} $\implies$ \eqref{eq:faktoryzacja1}]
We use the induction over $n$. We let us split the sum in equation
\eqref{eq:happynewyear} into the sum of terms with $\pi_{(1)}$ is
the maximal partition with only one block and the sum of all the
other terms; therefore
\begin{multline*}
\label{eq:happynewyear} k(\sigma_{1},\dots,\sigma_{n}) \E \big(
\Sigma_{\underbrace{1,\dots,1}_{l_1 \text{ times}} }\bullet\cdots
\bullet \Sigma_{\underbrace{1,\dots,1}_{l_n \text{
times}} } \big)= k^{\bullet}(\Sigma_{l_1},\dots,\Sigma_{l_n})- \\
\sum_{\pi_{(1)}, \pi_{(2)}}
k_{\pi_{(1)}}(\sigma_{1},\dots,\sigma_{n}) \ k^{\bullet}_{\pi_{(2)}}
(\Sigma_{\underbrace{1,\dots,1}_{l_1 \text{ times}}
},\dots,\Sigma_{\underbrace{1,\dots,1}_{l_n \text{ times}} }),
\end{multline*}
where the sum runs over partitions $\pi_{(1)},\pi_{(2)}$ such as in
\eqref{eq:happynewyear} with the additional constraint that
$\pi_{(1)}$ is not equal to the maximal partition with only one
block. We use the inductive hypothesis and estimate the summands on
the right-hand side in the same way as in the proof of the opposite
implication above.
\end{proof}

\subsection{Proof of the main theorem: equivalence \eqref{eq:faktoryzacja3} $\iff$
\eqref{eq:faktoryzacja4}}

\begin{proof}[Proof of equivalence \eqref{eq:faktoryzacja3} $\iff$
\eqref{eq:faktoryzacja4}] It is enough to apply Corollary
\ref{cor:roznegradacjejednadegradacja}.
\end{proof}

The remaining part of Theorem \ref{theo:main2}, namely
$$ \lim_{q\to\infty} \Cov(R_{l_1+1},R_{l_2+1} ) q^{-\frac{l_1+l_2}{2}}
=\lim_{q\to\infty}\Cov(\Sigma_{l_1},\Sigma_{l_2})
q^{-\frac{l_1+l_2}{2} } $$ follows from the relation
\eqref{eq:wolnekumulantyiklasy} between free cumulants and conjugacy
classes.

\section{Proof of Theorem \ref{theo:kumulantykombinatoryczne}}
\label{sec:Proof-Theorem-?????} This section is devoted to the proof
of Theorem \ref{theo:kumulantykombinatoryczne} which is the only
missing component in the proof of the main theorem. We will use some
tools presented in our recent work \cite{Sniady2003pushing}.

\subsection{Cumulants}
If $\mu=\{\mu_1,\dots,\mu_k\}$ is a partition of the set
$\{1,\dots,n\}$ with blocks $\mu_i=\{\mu_{i,1},\dots,\mu_{i,m(i)}
\}$ we define partition--indexed cumulants given by a multiplicative
extension of the usual cumulants:
\begin{equation}
\label{eq:multiplikatywnie} k_{\mu}(X_1,\dots,X_n)= \prod_i
k(X_{\mu_{i,1} }, X_{\mu_{i,2}},\dots, X_{\mu_{i,m(i)}});
\end{equation}
for example
$$  k_{ \{1,3,4\},\{2,5\}} (X_1,X_2,X_3,X_4,X_5)= k (X_1,X_3,X_4)\
k (X_2,X_5). $$ In this article we will use the following property
of cumulants: it turns out that the cumulants are implicitly
determined by a sequence of relations
\begin{equation}
\label{eq:moment-cumulant} \E ( X_1 \cdots X_n )= \sum_{\mu}
k_{\mu}(X_1,\dots,X_n),
\end{equation}
where the sum runs over all partitions $\mu$ of the set
$\{1,\dots,n\}$.

\subsection{Partition--indexed conjugacy classes}
In the following we present some constructions on partitions of the
set $X=\{1,2,\dots,N\}$. However, it should be understood that by a
change of labels these constructions can be performed for any finite
ordered set $X$.

We consider a matrix $J$, the entries of which belong to $\C(S_q)$,
the symmetric group algebra:
$$ J=\left[
\begin{matrix}   0 & (1,2)
& \dots &  (1,q) & 1 \\
 (2,1) & 0 &
 \dots & (2,q) & 1 \\
 \vdots & \vdots &
 \ddots & \vdots & \vdots \\
 (q,1) & (q,2)
 & \dots & 0 & 1 \\
 1 & 1 & 
  \dots  & 1 & 0
 \end{matrix} \right]\in \M_{q+1}(\C)\otimes \C(S_{q}). $$
Except for the last row, the last column and the diagonal, the entry
in the $i$-th row and the $j$-th column is equal to the
transposition interchanging $i$ and $j$.

Let $p=(p_1,\dots,p_l)$ be a sequence with
$p_1,\dots,p_l\in\{1,\dots,q+1\}$ and let $\pi$ be a partition of
the set $\{1,\dots,l\}$. We say that $p\sim\pi$ if for any $1\leq
i,j\leq l$ the equality $p_i=p_j$ holds if and only if $i$ and $j$
belong to the same block of the partition $\pi$. We define
\begin{equation}
\label{equ:noweklasy2} \Sigma_{\pi}=\sum_{\substack{p\sim\pi\\
p_l=q+1}} J_{p_1 p_2} J_{p_2 p_3} \cdots J_{p_{l-1} p_l} J_{p_l
p_1}\in \C(S_q).
\end{equation}
We will treat each summand as a partial permutation with a support
$\{p_1,\dots,p_l\}\setminus\{q+1\}$. Some partial results concerning
expressions of this form were obtained by Biane \cite{Biane1998}. We
can show \cite{Sniady2003pushing} that
$\Sigma_\pi=\Sigma_{k_1,\dots,k_t}$ for some integers
$k_1,\dots,k_t\geq 1$ which will be presented explicitly in Section
\ref{subsec:explicit}. For this reason we call $\Sigma_{\pi}$ a
partition-indexed conjugacy class.

\subsection{Products of conjugacy classes}
\begin{theorem}
\label{claim:mnozenie2} \label{claim:mnozenie} Let
$i_1<\cdots<i_{n+1}$ be integers and for each $1\leq s\leq n$ let
$\pi_s$ be a partition of the set
$\rho_s=\{i_{s}+1,i_{s}+2,\dots,i_{s+1}\}$. We denote
$$ \pi_1 \bullet \cdots \bullet \pi_n=   (\pi_1 \cup \cdots \cup \pi_n) \vee
\big\{ \{i_2,i_3,\dots,i_{n+1} \} \big\}. $$ Then
\begin{equation}
\label{eq:formulanailoczynrozlaczny} \Sigma_{\pi_1} \bullet \cdots
\bullet \Sigma_{\pi_n}= \Sigma_{(\pi_1 \bullet \cdots \bullet
\pi_n)}.
\end{equation}

Furthermore,
\begin{equation}
\label{eq:formulanailoczyn} \Sigma_{\pi_1} \cdots
\Sigma_{\pi_n}=\sum_{\sigma} \Sigma_{\sigma},
\end{equation}
where the sum on the right-hand side runs over all partitions
$\sigma$ of the set $\{i_1+1,i_1+2,\dots,i_{n+1}\}$ such that
\begin{enumerate}
 \item for any $a,b\in\rho_s$, $1\leq s\leq n$ we have that $a$ and
$b$ are connected by $\sigma$ if and only if they are connected by
$\pi_s$,
 \item elements $i_2,i_3,\dots,i_{n+1}$ belong to the same block of
 $\sigma$.
\end{enumerate}

Furthermore,
\begin{equation}
\label{eq:kumulantyklas} k^{\text{\rm
id}}(\Sigma_{\pi_1},\dots,\Sigma_{\pi_n})=\sum_{\sigma}
\Sigma_{\sigma},
\end{equation}
where the sum on the right-hand side runs over all partitions
$\sigma$ of the set $\{i_1+1,i_1+2,\dots,i_{n+1}\}$ such that
\begin{enumerate}
 \item for any $a,b\in\rho_s$, $1\leq s\leq n$ we have that $a$ and
$b$ are connected by $\sigma$ if and only if they are connected by
$\pi_s$,
 \item elements $i_2,i_3,\dots,i_{n+1}$ belong to the same block of
 $\sigma$; we will denote this block by $\sigma_1$,
 \item $(\sigma\setminus\sigma_1) \vee \{\rho_1,\dots,\rho_n\}$ is the
 maximal partition with only one block.
\end{enumerate}

\end{theorem}
\begin{proof}
Equations \eqref{eq:formulanailoczynrozlaczny} and
\eqref{eq:formulanailoczyn} follow immediately from the definition
\eqref{equ:noweklasy2}.

We shall treat \eqref{eq:kumulantyklas} as a definition of the
left-hand side and we shall verify that so defined cumulants fulfill
the defining relation of cumulants \eqref{eq:moment-cumulant}.
Before we do this we need to compute the corresponding
partition-indexed cumulants. Let $\mu$ be a partition of
$\{1,\dots,n\}$; we denote by $\tilde{\mu}$ a partition of the set
$\{i_1+1,i_1+2,\dots,i_{n+1}\}$ such that $a\in \rho_s$ and
$b\in\rho_t$ belong to the same block of $\tilde{\mu}$ if and only
if $s$ and $t$ belong to the same block of $\mu$. Then
\eqref{eq:kumulantyklas} implies
\begin{equation}
\label{eq:kumulantymultiplikatywnie} k_{\mu}^{\text{\rm
id}}(\Sigma_{\pi_1},\dots,\Sigma_{\pi_n})=\sum_{\sigma}
\Sigma_{\sigma},
\end{equation}
where the sum on the right-hand side runs over all partitions
$\sigma$ of the set $\{i_1+1,i_1+2,\dots,i_{n+1}\}$ such that
\begin{enumerate}
 \item for any $a,b\in\rho_s$, $1\leq s\leq n$ we have that $a$ and
$b$ are connected by $\sigma$ if and only if they are connected by
$\pi_s$,
 \item elements $i_2,i_3,\dots,i_{n+1}$ belong to the same block of
 $\sigma$; we will denote this block by $\sigma_1$,
 \item $(\sigma\setminus\sigma_1) \vee
 \{\rho_1,\dots,\rho_n\}=\tilde{\mu}$.
\end{enumerate}
Equation \eqref{eq:kumulantymultiplikatywnie} implies immediately
\eqref{eq:moment-cumulant}.
\end{proof}

For example, for
\begin{multline}
\label{eq:przykladmnozenie} \rho_1=\{1,2,3,4\},\
\rho_2=\{5,6,7,8\},\\  \pi_1=\big\{ \{1,3\},\{2,4\} \big\},
\pi_2=\big\{ \{5,7\}, \{6\},\{8\} \big\},\
\end{multline}
the above theorem states that
$$\Sigma_{\pi_1}  \Sigma_{\pi_2} =
\Sigma_{\big\{ \{1,3\},\{2,4,8\},\{5,7\},\{6\} \big\}}+
\Sigma_{\big\{ \{1,3,6\},\{2,4,8\},\{5,7\} \big\}}+  \Sigma_{\big\{
\{1,3,5,7\},\{2,4,8\},\{6\} \big\}}$$ and
$$k^{\text{\rm
id}}(\Sigma_{\pi_1},  \Sigma_{\pi_2}) = \Sigma_{\big\{
\{1,3,6\},\{2,4,8\},\{5,7\} \big\}}+ \Sigma_{\big\{
\{1,3,5,7\},\{2,4,8\},\{6\} \big\}};$$ the readers acquainted with
the results of Section \ref{subsec:explicit} may check that it is
equivalent to
$$ \Sigma_{1} \Sigma_{1,1} = \Sigma_{1,1,1}+ \Sigma_{1,1} +
\Sigma_{1,1} $$ and
$$ k^{\text{\rm
id}}(\Sigma_{1}, \Sigma_{1,1}) = \Sigma_{1,1} + \Sigma_{1,1}.
$$

\subsection{Geometric interpretation of the degree of $\Sigma_{\pi}$}
\label{subsec:zaklejanie}

\begin{figure}[tb]
\includegraphics{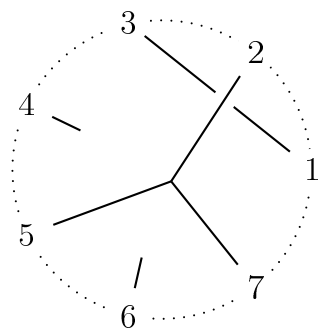}
\caption[A graphical representation of a partition]{Graphical
representation of a partition $\big\{ \{1,3\},\{2,5,7\},\{4\},\{6\}
\big\}$.} \label{fig:przecinajacapartycja}
\end{figure}

It is very useful to represent partitions graphically by arranging
the elements of the set $X=\{1,\dots,n\}$ counterclockwise on a
circle and joining elements of the same block by a line, as it can
be seen on Figure \ref{fig:przecinajacapartycja}.

\begin{figure}[tb]
\includegraphics{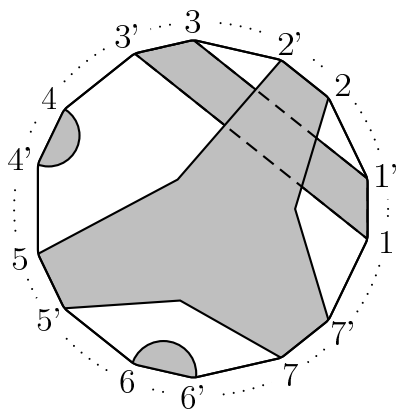}
\caption[The first collection of discs for $\pi$ from Figure
\ref{fig:przecinajacapartycja}.]{The first collection of discs for
partition $\pi$ from Figure \ref{fig:przecinajacapartycja}.}
\label{fig:firstcollection}
\end{figure}

We consider a large sphere with a small circular hole. The boundary
of this hole is the circle mentioned above. Let us draw the blocks
of the partition $\pi$ with a fat pen; in this way each block
becomes a disc glued to the boundary of the hole, cf Figure
\ref{fig:firstcollection}.

After gluing the first collection of discs, our sphere becomes a
surface with a number of holes. The boundary of each hole is a
circle and we shall glue this hole with a disc from the second
collection. Thus we obtained an orientable surface without a
boundary. We call the genus of this surface the genus of the
partition $\pi$ and denote it by $\genus_{\pi}$.


The following result was proved in our previous work
\cite{Sniady2003pushing}.
\begin{proposition}
\label{prop:genusowatosc} For any partition $\pi$ of an $n$--element
set
\begin{equation}
\label{eq:genusowatosc} \deg \Sigma_{\pi}=n-2 \genus_\pi.
\end{equation}
\end{proposition}


\subsection{Geometric interpretation of Theorem \ref{claim:mnozenie}}
\label{subsec:geometricmultiplication} 
\begin{figure}[bt]
\includegraphics{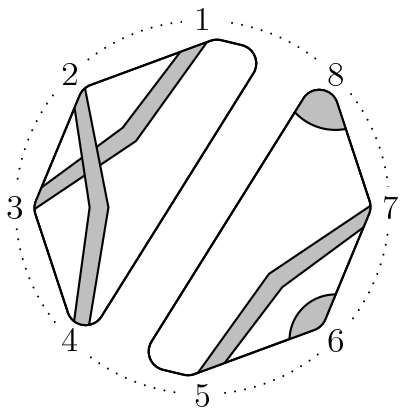}
\caption{Graphical representation of example
\eqref{eq:przykladmnozenie}.} \label{fig:mnozeniepartycji}
\end{figure}

We will use the notations of
Theorem \ref{claim:mnozenie}. On the surface of a large sphere we
draw a small circle on which we mark counterclockwise points
$i_1+1,i_1+2,\dots,i_{r+1}$.
Inside the circle we cut $r$ holes; for any $1\leq s\leq r$ the
corresponding hole has a shape of a disc, the boundary of which
passes through the points from the block $\rho_s$. For every $1\leq
s\leq r$ the partition $\pi_s$ connects some points on the boundary
of the hole $\rho_s$ and this situation corresponds exactly to the
case we considered in Section \ref{subsec:zaklejanie}. We shall glue
to the hole $\rho_s$ only the first collection of discs that we
considered in Section \ref{subsec:zaklejanie}, i.e.~the discs which
correspond to the blocks of the partition $\pi_s$. Thus we obtained
a number of holes with a collection of glued discs (cf Figure
\ref{fig:mnozeniepartycji}).

\begin{figure}[bt]
\includegraphics{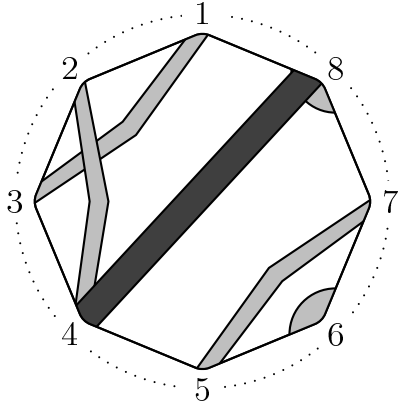}
\caption{Figure \ref{fig:mnozeniepartycji} after inflating small
holes.} \label{fig:mnozeniepartycji2}
\end{figure}

When we inflate the original small holes inside the circle we may
think about this picture alternatively: instead of $r$ small holes
we have a big one (in the shape of the circle) but some arcs on its
boundary are glued by an additional disc (on Figure
\ref{fig:mnozeniepartycji2} drawn in black) glued to vertices
$i_2,i_3,\dots,i_{r+1}$. Furthermore we still have a collection of
all discs (on Figure \ref{fig:mnozeniepartycji2} drawn in gray)
corresponding to partitions $\pi_s$. We merge the additional black
disc to a gray disc from this collection if they touch the same
vertex. After this merging the collection of discs corresponds to
the partition $\pi_1 \bullet \cdots \bullet \pi_r=\big\{
\{i_2,i_3,\dots,i_{r+1} \} \big\} \vee (\pi_1 \cup \cdots \cup
\pi_r)$ which appears in the formula
\eqref{eq:formulanailoczynrozlaczny} for the disjoint product.

The last step is to consider all ways of merging of the discs (or
equivalently: all partitions $\sigma\geq \pi_1 \bullet \cdots
\bullet \pi_r$) with the property that any two vertices that were
lying on the boundary of the same small hole $\rho_s$ if were not
connected by a disc from the collection $\pi_s$ then they also
cannot be connected after all mergings. In this way we obtain all
partitions which contribute to \eqref{eq:formulanailoczyn}.

By splitting the holes (we recall that each hole corresponds to some
set $\rho_i$) into some new holes we can view the surface associated
to the partition $\rho$ as a sphere with a number of new holes glued
in pairs by handles and the number of these handles is equal to the
genus of the surface. In this way we obtain a graph $G_\sigma$ the
vertices of which correspond to the old holes (or, equivalently,
sets $\rho_i$) and the edges correspond to the handles between new
holes. Of course the above construction can be sometimes performed
in many different ways but we do not mind it. Note that multiple
connections between vertices are allowed. Also, a vertex can be
connected with itself and the number of such loops is equal to the
genus of the corresponding partitions $\pi_i$. The genus of $\sigma$
is equal to the total number of edges in $G_{\sigma}$.

In order to obtain all partitions $\sigma$ which contribute to
\eqref{eq:kumulantyklas} we should restrict our attention only to
$\sigma$ such that the graph $G_{\sigma}$ is connected. A connected
graph with $n$ vertices has at least $n-1$ non-loop edges therefore
\begin{equation}
\label{eq:nierownosc} \genus_{\sigma}\geq
\genus_{\pi_1}+\cdots+\genus_{\pi_n}+(n-1).
\end{equation}

\begin{proof}[Proof of Theorem \ref{theo:kumulantykombinatoryczne}]
In order to prove \eqref{eq:kumulantykombinatoryczne} for the
special case when $x_i=\Sigma_{\pi_i}$ it is enough to apply
Proposition \ref{prop:genusowatosc} and \eqref{eq:nierownosc}.

It remains now to prove \eqref{eq:a-tak-dokladnie?} and we shall do
it in the following. We use the notations of Theorem
\ref{claim:mnozenie2} if $n=2$ and $\pi_1,\pi_2$ are trivial
partitions (every blocks consists of a single element). Let $\sigma$
be a partition which contributes to \eqref{eq:kumulantyklas} with
the minimal possible genus, namely $\genus_{\sigma}=1$. It follows
that one of the blocks of $\sigma$ is equal to $\{i_2,i_3\}$.
Secondly, some of the elements of the set $\{i_1+1,\dots,i_2-1\}$
(but at least one of them) are paired with some of the elements of
the set $\{i_2+1,\dots,i_1+i_3-1\}$; this pairing however is not
arbitrary. Let us travel counterclockwise along the boundary of the
hole corresponding to the block $\pi_1$; in other words we visit the
vertices in the order $i_1+1\rightarrow i_1+
2\rightarrow\cdots\rightarrow i_2-1\rightarrow i_1+1$. Some of these
vertices are paired by $\sigma$ with some of the elements of
$\{i_2+1,\dots,i_3-1\}$ corresponding to the other hole; let us have
a look in which order these counterparts appear during our walk.
From the very definition of the genus of a partition we know that it
is possible to draw the blocks of $\sigma\setminus \big\{
\{i_2,i_3\} \big\}$ on the surface of the handle in such a way that
lines do not cross.  It follows that these counterparts will be
visited in the clockwise order $i_3-1\rightarrow i_3-2\rightarrow
\cdots\rightarrow i_2+1 \rightarrow i_3-1$, cf Figure
\ref{fig:cake2}.

Let $a_1,\dots,a_r$ denote the distances (counted cyclically
counterclockwise) between consecutive elements of the pairs in the
first hole and let $b_1,\dots,b_r$ denote the distances (counted
cyclically counterclockwise) between consecutive elements of the
pairs in the second hole. Readers acquainted with the results of
Section \ref{subsec:explicit} will see that
$$ \Sigma_{\sigma}= \Sigma_{a_1+b_1-1,\dots,a_r+b_r-1}.$$
It is easy to see that $a_1+\cdots+a_r=i_2-i_1$ and
$b_1+\cdots+b_r=i_3-i_2$. Such sequences $(a_i)$ and $(b_i)$
uniquely determine $\sigma$ once we specify the first element in
each cycle (there are $i_2-i_1-1$ choices for the first one and
$i_3-i_2-1$ choices for the second one). Partition $\sigma$ can be
represented like this in $r$ different ways which correspond to the
cyclic rotations of the sequences $(a_i)$ and $(b_i)$. Since
$\Sigma_{\pi_1}=\Sigma_{i_2-i_1-1}$ and
$\Sigma_{\pi_2}=\Sigma_{i_3-i_2-1}$ this finishes the proof of
\eqref{eq:a-tak-dokladnie?}.
\begin{figure}[tb]
\includegraphics{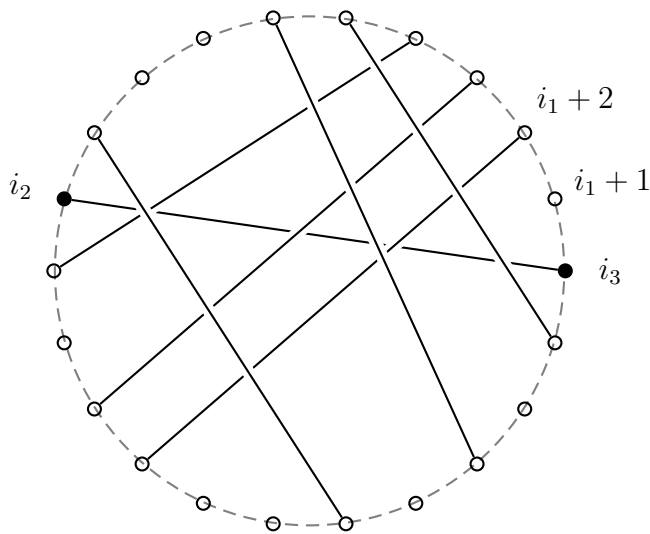}
\caption{Partitions $\sigma$ with genus $1$ which contribute to
\eqref{eq:kumulantyklas}.} \label{fig:cake2}
\end{figure}
\end{proof}

%
%

\subsection{Fat partitions} \index{fat partition}
\index{partition!fat} \index{$Pi$@$\pi_{\fat}$} \label{subsec:fat}
Let $\pi=\{\pi_1,\dots,\pi_r\}$ be a partition of the set
$\{1,\dots,n\}$. For every $1\leq s\leq r$ let
$\pi_s=\{\pi_{s,1},\dots,\pi_{s,l_s}\}$ with
$\pi_{s,1}<\cdots<\pi_{s,l_s}$. We define $\pi_{\fat}$, called fat
partition of $\pi$, to be a pair partition of the $2n$--element
ordered set $\{1,1',2,2',\dots,n,n'\}$ given by
$$\pi_{\fat}=
\big\{  \{ \pi_{s,t}',\pi_{s,t+1}  \} :  1\leq s \leq r \text{ and }
1 \leq t \leq l_s  \big\}, $$ where it should be understood that
$\pi_{s,l_s+1}=\pi_{s,1}$.

\begin{figure}[tb]
\includegraphics{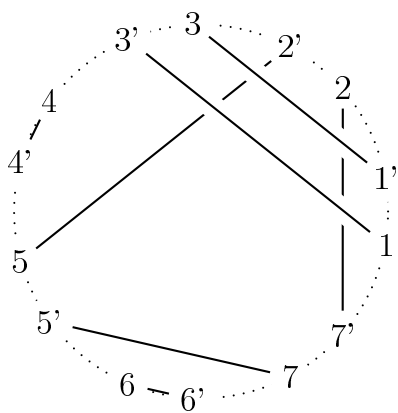}
\caption[Fat partition $\pi_{\fat}$]{The fat partition $\pi_{\fat}$
corresponding to the partition $\pi$ from Figure\
\ref{fig:przecinajacapartycja}.} \label{fig:tlustapartycja}
\end{figure}

This operation can be easily described graphically as follows: we
draw the blocks of the partition with a fat pen and take the
boundary of each block, as it can be seen on Figure
\ref{fig:tlustapartycja}. This boundary is a collection of lines
hence it is a pair partition. However, every vertex
$k\in\{1,\dots,n\}$ of the original partition $\pi$ has to be
replaced by its `right' and `left' copy (denoted respectively by $k$
and $k'$).



\subsection{Explicit form of the partition--indexed conjugacy class indicator $\Sigma_\pi$}
\label{subsec:explicit} Let $\pi$ be a partition of the set
$\{1,\dots,n\}$. Since the fat partition $\pi_{\fat}$ connects every
element of the set $\{1',2',\dots,n'\}$ with exactly one element of
the set $\{1,2,\dots,n\}$, we can view $\pi_{\fat}$ as a bijection
$\pi_{\fat}:\{1',2',\dots,n'\}\rightarrow\{1,2,\dots,n\}$. We also
consider a bijection
$c:\{1,2,\dots,n\}\rightarrow\{1',2',\dots,n'\}$ given by
$\dots,3\mapsto 2', 2\mapsto 1', 1\mapsto n', n\mapsto
(n-1)',\dots$. Finally, we consider a permutation $\pi_{\fat}\circ
c$ of the set $\{1,2,\dots,n\}$.

\begin{figure}[bt]
\includegraphics{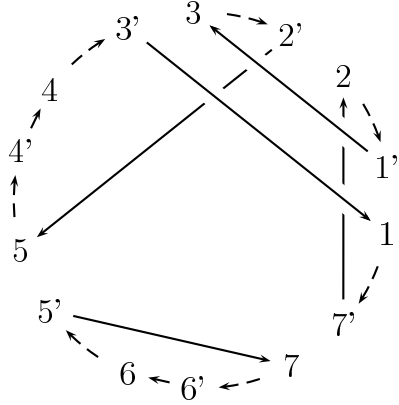}
\caption[Bijection corresponding to $\pi_{\fat}$ from Figure
\ref{fig:tlustapartycja}]{Bijection corresponding to the partition
$\pi_{\fat}$ from Figure \ref{fig:tlustapartycja} plotted with a
solid line and the bijection $c$ plotted with a dashed line. }
\label{fig:tlustapartycja2}
\end{figure}

For example, for the partition $\pi$ given by Figure
\ref{fig:przecinajacapartycja} the composition $\pi_{\fat}\circ c$
has a cycle decomposition $(1,2,3,5,{4})({6},7)$, as it can be seen
from Figure \ref{fig:tlustapartycja2}.

\begin{figure}[bt]
\includegraphics{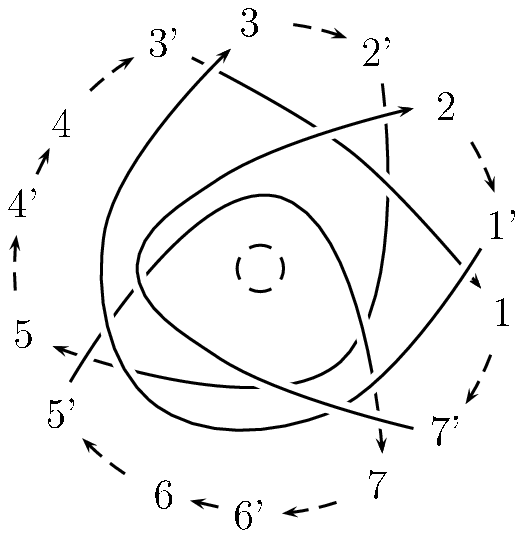}
\caption[Version of Figure \ref{fig:tlustapartycja2} in which all
lines wind clockwise]{A version of Figure \ref{fig:tlustapartycja2}
in which all lines wind clockwise around the central disc.}
\label{fig:tlustapartycja3}
\end{figure}

We decompose the permutation
$$
\pi_{\fat}\circ c=(b_{1,1},b_{1,2},\dots,b_{1,j_1}) \cdots
(b_{t,1},\dots,b_{t,j_t})
$$
as a product of disjoint cycles. Every cycle
$b_s=(b_{s,1},\dots,b_{s,j_s})$ can be viewed as a closed clockwise
path on a circle and therefore one can compute how many times it
winds around the circle, cf Figure \ref{fig:tlustapartycja3}.

To a cycle $b_s$ we assign the number
\begin{multline*}
k_s=(\text{number of elements in a cycle } b_s)-\\ (\text{number of
clockwise winds of }b_s).
\end{multline*}

In the above example we have $b_1=(1,2,3,5,4)$, $b_2=(6,7)$ and
$k_1=2$, $k_2=1$, as it can be seen from Figure
\ref{fig:tlustapartycja3}, where all lines clockwise wind around the
central disc.

In our recent work \cite{Sniady2003pushing} we proved that
$$\Sigma_\pi=\Sigma_{k_1,\dots,k_t},$$
where $\Sigma_{k_1,\dots,k_t}$ on the right--hand side should be
understood as in Section \ref{subsec:definicjasigma}.

\section{Acknowledgments}

I thank Akihito Hora and Philippe Biane for many discussions and
remarks.

Research supported by State Committee for Scientific Research
(Komitet Bada\'n Naukowych) grant No.\ 2P03A00723; by EU Network
``QP-applications", contract HPRN-CT-2002-00729; by KBN-DAAD project
36/2003/2004.

\bibliographystyle{alpha}
\bibliography{biblio}

\end{document}